\documentclass[preprint]{imsart}

\RequirePackage[OT1]{fontenc}
\RequirePackage{amsthm,amsmath}
\RequirePackage[numbers]{natbib}
\RequirePackage[colorlinks,citecolor=blue,urlcolor=blue]{hyperref}

\RequirePackage{mathptmx}
\RequirePackage{amssymb}
\RequirePackage{mathrsfs}
\RequirePackage{amsfonts,dsfont}
\RequirePackage{graphicx}
\RequirePackage{psfrag}
\RequirePackage{booktabs}
\RequirePackage{url}
\RequirePackage[usenames]{color}

\newcommand{\R}{\mathbb{R}}
\newcommand{\N}{\mathbb{N}}
\newcommand{\E}{\mathbb{E}}

\newcommand{\dd}{\text{d}}
\newcommand{\1}{\mathds{1}}
\newcommand{\F}{\mathcal{F}}

\startlocaldefs
\numberwithin{equation}{section}
\usepackage[top=1.2in, bottom=1.2in, left=1.4in, right=1.4in]{geometry}
\theoremstyle{plain}
\newtheorem{thm}{Theorem}[section]

\newtheorem{lem}[thm]{Lemma}

\newtheorem{cor}[thm]{Corollary}

\newtheorem{assumption}[thm]{Assumption}
\endlocaldefs

\begin{document}

\begin{frontmatter}
\title{An efficient explicit full-discrete scheme for strong approximation of stochastic Allen-Cahn equation\thanksref{T1}}
\runtitle{Full discrete scheme for stochastic Allen-Cahn equation}
\thankstext{T1}{This work was supported by NSF of China (11671405, 11971488, 91630312), Natural
Science Foundation of Hunan Province (2020JJ2040, 2018JJ3628).
The author wants to thank the Tianyuan Mathematical Center in Northeast China
for the hospitality 
when this work was presented in a conference in June of 2018, hosted by the center.
The author  also  wants to thank Yuying Zhao for her kind help with excellent typesetting
and Ruisheng Qi and Meng Cai for their useful comments based on carefully reading the manuscript.
}

\author{\fnms{Xiaojie} \snm{Wang}
\ead[label=e1]{x.j.wang7@csu.edu.cn
}}
%
%
\address{School of Mathematics and Statistics,\\
Central South University, Changsha, Hunan, China\\
\printead{e1}
\\
\phantom{E-mail:\ }
}
\runauthor{X. Wang}



\begin{abstract}
In [Becker and Jentzen, Stoch. Proc. Appl. 129:28--69, 2019] and [Becker, et al., arXiv preprint arXiv:1711.02423, 2017],  
an explicit temporal semi-discretization scheme and a space-time full-discretization scheme were, respectively, introduced and analyzed for 
the additive noise-driven stochastic Allen-Cahn type equations, with strong convergence rates recovered.
The present work aims to propose a different explicit full-discrete scheme to numerically solve the stochastic Allen-Cahn equation 
with cubic nonlinearity, perturbed by additive space-time white noise. The approximation is easily implementable, performing the spatial 
discretization by a spectral Galerkin method and the temporal discretization by a kind of nonlinearity-tamed accelerated exponential integrator 
scheme. Error bounds in a strong sense are analyzed for both the spatial semi-discretization and the spatio-temporal full discretization, 
with convergence rates in both space and time explicitly identified. It turns out that the obtained convergence rate of 
the new scheme is, in the temporal direction, twice as high as existing ones in the literature. 
Numerical results are finally reported to confirm the previous theoretical findings.
\end{abstract}

\begin{keyword}[class=MSC]
65C30, 60H35,  60H15
\end{keyword}

\begin{keyword}
\kwd{stochastic Allen-Cahn equation}
\kwd{cubic nonlinearity}
\kwd{spectral Galerkin method}
\kwd{tamed exponential integrator scheme}
\kwd{strong convergence rate}
\end{keyword}

\end{frontmatter}

\section{Introduction}\label{sec:introduction}

As an active area of research, numerical study of evolutionary stochastic partial differential equations (SPDEs) 
has attracted increasing attention in the past decades 
(see, e.g.,  monographs \cite{kruse2014strong,lord2014introduction,jentzen2011taylor} and references therein). 
Albeit much progress has been made,  it is still not well-understood, especially for numerical analysis of 
SPDEs with non-globally Lipschitz nonlinearities. 
The present work attempts to make a contribution in this direction and examine
a space-time full-discretization scheme for a typical example of parabolic SPDEs with super-linearly growing nonlinearity, 
i.e., the stochastic Allen-Cahn equation.
%
The driving noise is a space-time white one, which is of special interest as it can best model the fluctuations 
generated by microscopic effects in a homogeneous physical system, including, for example, 
molecular collisions in gases and liquids, electric fluctuations in resistors \cite{gardiner1986handbook}.
A lot of researchers carried out numerical analysis of SPDEs subject to such noise, e.g., 
\cite{liu2003convergence, 
faris1982large,          
GI98, GI99,      
davie2001convergence,  
jentzen2009overcoming,
becker2016strong,          
blomker2013galerkin,     
hutzenthaler2016strong,      
Anton2017fully,          
DP09,DA10,                
jentzen2011higher,
brehier2018analysis,  
liu2018strong,            
cao2017approximating,
printems2001discretization,  
wang2014higher},      
to just mention a few.

Numerically solving the continuous problem on a computer forces us to perform 
both spatial and temporal discretizations. In space, we discretize SPDE \eqref{eq:SGL-concrete} 
by a spectral Galerkin method,
resulting in a system of finite dimensional stochastic differential equations (SDEs). Based on the spatial discretization,  
we propose a nonlinearity-tamed accelerated exponential time-stepping scheme given by \eqref{eq:full.Tamed-AEE}. 
The resulting approximation errors of both the spatial discretization and the space-time fully discrete scheme are carefully analyzed, 
with strong convergence rates successfully recovered.
More accurately, by $X(t_m)$ we denote the unique mild solution of the underlying  SPDE  taking values at temporal grid points 
$t_m = m \tau, m \in \{ 0, 1,..., M \}$ with uniform time step-size $\tau = \tfrac{T}{M} > 0$ and by $Y^{M,N}_{t_m}$ the numerical approximations 
of $X(t_m)$,  produced by the proposed fully discrete scheme. 
The approximation error measured in $L^p ( \Omega; H), p \in [2, \infty)$ reads
(cf. Theorem \ref{thm:full-discrete-scheme-error-bound}):
\begin{equation} \label{eq:Intro-numerical-main-result}
\sup_{ 0 \leq m \leq M}
\| X ( t_{m} ) -  Y^{M,N}_{t_m} \|_{L^p ( \Omega; H) }
     \leq
C 
\big ( N^{ - \beta } + \tau ^{ \beta } \big), 
\quad
\forall \, \beta \in (0, \tfrac12).
\end{equation}
Here $H := L^2 ((0, 1); \R )$ and the constant $C$ depends on $p, \beta, T$ and the initial value of the SPDE,
but does not depend on the discretization parameters $M, N$. 
%

Over the last two years, several research works were reported on numerical analysis of space-time white noise driven SPDEs with 
cubic (polynomial) nonlinearity \cite{becker2016strong,Becker2017strong,liu2018strong,brehier2018analysis,brehier2018strong}.
Becker and Jentzen \cite{becker2016strong} in 2016 introduced two nonlinearity-truncated Euler-type approximations
for pure time discretizations of stochastic Ginzburg Landau type equations 
with slightly more general polynomials. There a strong convergence rate of order almost $\frac{1}{4} $ is identified. 
More recently when the first preprint of this work was almost finished, we were aware of four other preprints 
\cite{Becker2017strong,brehier2018analysis,liu2018strong,brehier2018strong} submitted to arXiv, 
concerning with numerical approximations of 
similar SPDEs.
Becker, Gess, Jentzen and Kloeden \cite{Becker2017strong} 
proposed new types of truncated exponential Euler space-time fully discrete schemes for the same problem 
as in  \cite{becker2016strong} and derived strong convergence rates of order almost $\frac{1}{2} $ in space  
and order almost $\frac{1}{4} $ in time. 
Later, Br{\'e}hier and Gouden\`ege \cite{brehier2018analysis} and Br{\'e}hier, Cui and Hong \cite{brehier2018strong} analyzed 
some splitting time discretization schemes and obtained strong convergence rates of order $\frac14$. 
Liu and Qiao \cite{liu2018strong} investigated a spectral Galerkin backward implicit Euler full discretization,
with strong convergence rates of order almost $\frac{1}{2} $ in space  and order $\frac{1}{4} $ in time achieved.
As clearly implied by \eqref{eq:Intro-numerical-main-result}, the spatial convergence rate coincides with those in 
\cite{Becker2017strong,liu2018strong}, but the strong convergence rate of our time-stepping scheme 
can be of order almost $\frac12$, essentially twice as high as those in  
\cite{becker2016strong,Becker2017strong,liu2018strong,brehier2018analysis,brehier2018strong}. 
Despite getting involved with linear functionals of the noise process, 
the newly proposed scheme  is explicit, 
easily implementable and does not cost additional computational efforts
(see comments in section \ref{sec:numerical-result} for the implementation of the linear functionals of the noise process).

It is important to emphasize that, proving the error estimate \eqref{eq:Intro-numerical-main-result} is challenging,
confronted with two essential difficulties, one being to derive uniform a priori moment bounds for the numerical approximations
with super-linearly growing nonlinearity (see relevant comments in \cite{Becker2017strong}) 
and the other to recover the temporal convergence rate of order almost $\tfrac12$, 
instead of order (almost) $ \tfrac{1}{4} $ in the existing literature.  With regard to the former, we essentially rely on certain estimates for 
deterministic perturbed PDEs \eqref{eq:determ-perturbed-PDE}, as elaborated in subsection \ref{subsec:estimates-perturbed-PDE}.
The uniform a priori moment $L_\infty$-bounds are derived based on a certain bootstrap argument, by showing $ \E\big[ \1_{ \Omega_{ R^{ \tau },t_m} } 
\| Y^{M,N}_{t_m} \|_V^{p} \big] < \infty$ and $    \E
   \big [
   \1_{\Omega^c_{R^{\tau},t_m}} \|Y_{t_m}^{M,N} \|_V^p
   \big ] < \infty $,
   $ V : = C ( (0, 1), \R )$,
for  subevents $ \Omega_{ R^{ \tau },t_m} $  with $ R^{ \tau } $ depending on $\tau$ carefully chosen 
(see subsection \ref{subsec:a-priori-moment-numerical}).
The latter difficulty is addressed by fully exploiting
the global monotonicity condition on the nonlinearity 
in conjunction with smoothing property of the analytic semigroup, commutativity properties of the nonlinearity and
improved temporal H\"{o}lder regularity results in negative Sobolev spaces (consult subsection \ref{subset:negative-sobleve}
and particularly the treatment of $J_1$ in the proof of Theorem \ref{thm:full-discrete-scheme-error-bound} for details).
%

Furthermore, we would like to point out that the improvement of convergence rate is essentially 
credited to fully preserving the stochastic convolution in the time-stepping scheme \eqref{eq:full.Tamed-AEE}.
Such a kind of accelerating technique was used by Jentzen and Kloeden \cite{jentzen2009overcoming},
to solve nearly linear parabolic SPDEs 
and has been further examined and extended in different settings \cite{jentzen2011efficient,wang2015note, 
wang2014higher,qi2017accelerated,lord2016modified}, where 
a globally Lipschitz condition imposed on nonlinearity is indispensable in the error analysis.
%
When the nonlinearity grows super-linearly and the globally Lipschitz condition is thus violated,  one can in general not expect
the usual accelerated exponential time-stepping schemes converge in the strong sense,
based on the observations in 
\cite{beccari2019strong, hutzenthaler2011strong}.  To address this issue,  
we introduce a taming technique previously used in \cite{hutzenthaler2012strong,wang2013tamed,Tretyakov2013fundamental,hutzenthaler15MEMAMS,hutzenthaler2013exponential} 
for ordinary SDEs, and propose a nonlinearity-tamed version of accelerated exponential Euler scheme for the time discretization.
Although the idea of constructing the explicit scheme in this paper is inspired by the aforementioned works, 
the approach of the error analysis for the full discretization of the stochastic Allen-Cahn equation with space-time white noise is original
and much more involved than that in both the finite dimensional non-globally Lipschitz SDE setting and the globally Lipschitz SPDE 
{\color{black}{setting}} (see section \ref{sect:full-discretization}).
Moreover, our approach is much easier than that in \cite{Becker2017strong}, 
also treating explicit full-discrete schemes for Allen-Cahn type SPDEs.
%

Finally, we mention that, just one spatial dimension is considered here because 
the space-time white noise driven SPDE only allows for a mild solution 
with a positive (but very low) order of regularity in one spatial dimension.
It is because of the low order of regularity that the error analysis becomes difficult.
Further, the error analysis of the space-time full discretization is significantly more involved than 
that of the pure time semi-discretization (compare \cite{becker2016strong} and \cite{Becker2017strong} 
and see comments therein).
Strong convergence analysis of numerical methods for smoother noise (e.g., trace-class noise) driven 
stochastic Allen-Cahn equation in multiple spatial dimensions has been done in our recent work \cite{qi2019optimal}
and further strong and weak approximations of parabolic SPDEs with non-globally Lipschitz nonlinearity 
will be our forthcoming works
(see also, e.g.,
%
\cite{brehier2018weak,
sauer2015lattice,      
kovacs15backward, kovacs15discretisation,   
majee2017optimal,     
feng2017finite,
gyongy2016convergence,   
jentzen2015strong,          
hutzenthaler2014perturbation
}
for revalent topics).

The rest of this paper is organized as follows. In the next section we collect some basic facts and
present the well-posedness of the stochastic problem under given assumptions. 
Section \ref{sect:spatial-discret} and Section \ref{sect:full-discretization} are, respectively, 
devoted to the analysis of strong convergence rates for both the spatial semi-discretization and  
the spatio-temporal full discretization of the underlying SPDEs.   
Numerical results are included in section \ref{sec:numerical-result}  to test previous theoretical findings.

%


\section{Well-posedness of the stochastic problem}
\label{sec:well-posedness-regularity}
%
%
%
Throughout this article, we are interested in the additive space-time white noise driven 
stochastic Allen-Cahn equation with cubic nonlinearity, described by
\begin{equation}\label{eq:SGL-concrete}
\left\{
    \begin{array}{lll}
    \frac{\partial u }{\partial t }(t,x) = \frac{\partial^2 u }{\partial x^2 } (t,x) + f(u(t,x)) + \dot{W} (t, x),  \: x \in D,
     \ t \in (0, T],
    \\
     u(0, x) = u_0(x), \: x \in D,
     \\
     u(t, 0) = u (t, 1) = 0,  \: t \in (0, T].
    \end{array}\right.
\end{equation}
Here $D := (0, 1)$,  $ T > 0$, $f \colon \mathbb{R} \rightarrow \mathbb{R}$ is given by 
$
f (v) =  a_3 v^3  + a_2 v^2 + a_1 v + a_0, \,
a_{3} < 0, \, a_2, a_1, a_0, v \in \R,
$
and $\dot{W} (t, \cdot )$ stands for a formal time derivative of a cylindrical I-Wiener process \cite{da2014stochastic}.
%
%
In order to define a mild solution of \eqref{eq:SGL-concrete} following the semigroup approach in \cite{da2014stochastic},
we attempt to put everything into an abstract framework.
Given a real separable Hilbert space $(H, \langle \cdot, \cdot \rangle, \|\cdot\| )$ with $\|\cdot\| = \langle \cdot, \cdot \rangle^{\frac{1}{2}}$,
by $\mathcal{L}(H)$ we denote the space of bounded linear operators from 
$H$ to $H$ endowed with the usual operator norm $\| \cdot \|_{\mathcal{L}(H)}$.
Additionally, we denote by $\mathcal{L}_2(H) \subset \mathcal{L}(H)$ the subspace
consisting of all Hilbert-Schmidt operators from $H$ to $H$ \cite{da2014stochastic}. 
It is known that $\mathcal{L}_2(H)$ 
is a separable Hilbert space, equipped with the scalar product
$
\langle \Gamma_1, \Gamma_2  \rangle_{\mathcal{L}_2(H)} 
: = 
\sum_{n \in \N} \langle \Gamma_1 \eta_n, \Gamma_2 \eta_n  \rangle,
$
and norm
$
\| \Gamma \|_{ \mathcal{L}_2(H) } 
:=
\big (
\sum_{n \in \N} \| \Gamma \eta_n \|^2
\big )^{\frac12},
$
independent of the particular choice of orthogonal basis $\{\eta_n\}_{n \in \N}$ of $H$.
Below we sometimes write $\mathcal{L}_2 := \mathcal{L}_2(H)$ for brevity.
If $\Gamma \in \mathcal{L}(H)$ and $\Gamma_1, \Gamma_2 \in \mathcal{L}_2(H)$, then
$
| \langle \Gamma_1, \Gamma_2  \rangle_{\mathcal{L}_2(H)} | 
\leq 
\| \Gamma_1 \|_{\mathcal{L}_2(H)} 
\| \Gamma_2  \|_{\mathcal{L}_2(H)},
\| \Gamma \Gamma_1\|_{\mathcal{L}_2(H)} \leq \|\Gamma\|_{\mathcal{L}(H)} \|\Gamma_1\|_{\mathcal{L}_2(H)}.
$
By $L^{ \gamma } ( D; \R ), \gamma \geq 1$ ($L^{ \gamma } ( D )$ for short) 
we denote a Banach space consisting of $\gamma$-times integrable functions 
and by $ V : = C ( D, \R )$ a Banach space of continuous functions with usual norms.
%
%
%
%
%
%
%
%
%
To reformulate \eqref{eq:SGL-concrete} as an abstract problem,  we make the following assumptions.
\begin{assumption}[Linear operator $A$]\label{ass:A}
Denote $D : = (0, 1)$ and let $H = L^2 (D; \R )$  be a real separable Hilbert space, 
equipped with usual product  $\langle \cdot, \cdot \rangle$
and norm $\|\cdot \| = \langle \cdot, \cdot \rangle^{\frac12}$.
Let $ - A \colon \text{dom} (A) \subset H \rightarrow H$ be the Laplacian with homogeneous Dirichlet boundary conditions, 
defined by $ - A u = \Delta u$, $u \in \text{dom}(A) := H^2 \cap H_0^1$.
\end{assumption}
The above setting assures that there exists an increasing sequence of real numbers $\lambda_i=\pi^2 i^2, i \in \mathbb{N}$ 
and an orthonormal basis $\{e_i ( x ) =\sqrt{2}\sin(i \pi x), \,x
\in(0,1)\}_{i \in \mathbb{N}}$ such that $A e_i = \lambda_i e_i$.
In particular, the linear unbounded operator $A$ is positive, i.e., $ \langle - A v, v \rangle \leq - \lambda_1 \| v \|^2$,
for all $ v \in  \text{dom}(A)$.
Moreover, $-A$ generates an analytic semigroup $E(t) = e^{-t A}, t \geq 0$ on $H$
%
and we can define the fractional powers of $A$, i.e., $A^\gamma, \gamma \in \mathbb{R}$ 
and the Hilbert space $\dot{H}^{\gamma} :=\text{dom}(A^{\frac{\gamma}{2}})$, equipped with inner product  
$\langle \cdot, \cdot \rangle_\gamma : = \langle A^{\frac{\gamma}{2}} \cdot, 
A^{\frac{\gamma}{2}} \cdot \rangle$ 
and norm $\| \cdot \|_{\gamma} = \langle \cdot, \cdot \rangle_\gamma^{\frac12} $ 
\cite[Appendix B.2]{kruse2014strong}. Moreover, $\dot{H}^0 = H$ and 
$\dot{H}^{\gamma} \subset \dot{H}^{\delta},  \gamma \geq \delta$. It is well-known that \cite{pazy1983semigroups}, 
for a positive constant $c >0$,
\begin{equation}
\begin{split}
\label{eq:E.Inequality}
\| A^\gamma E(t)\|_{\mathcal{L}(H)} \leq& c t^{-\gamma}, \quad t >0, \gamma \geq 0,
\\
\| A^{-\rho} (I-E(t))\|_{\mathcal{L}(H)} \leq& c t^{\rho}, \quad t >0, \rho \in [0,1],
\end{split}
\end{equation}
%
and accordingly one can verify that \cite[Lemma B.9 (iii)]{kruse2014strong}
\begin{equation} \label{eq:semigroup-integral}
\int_s^t \| A^{\frac{\rho}{2}} E ( t - r ) u \|^2 \, \dd r 
\leq
c 
( t -  s )^{ 1 - \rho }
\| u \|^2,
\quad
\forall u \in H, 
\,
\rho \in [0, 1],
\,
0\leq s \leq t.
\end{equation}
Here the constant $c$, might depending on $\gamma$, can be chosen to be independent of $\rho \in [0, 1]$.
Also, it is evident to see that
\begin{equation}
\label{eq:AQ_condition}
\|A^{\frac{\beta-1}{2}}  \|_{\mathcal{L}_2 (H) }  
{\color{black}{
=
\pi^{\beta - 1} 
\Big(
\sum_{i \in \N} i^{ 2 ( \beta - 1) }
\Big)^\frac12
}}
< \infty, \quad  \text{ for any } \beta < \tfrac12.
\end{equation}

\begin{assumption}[Nonlinearity]\label{ass:F}
Let $F \colon L^{6} ( D; \R) \rightarrow H$ be  a deterministic mapping defined by
%
%
\begin{equation*}
F (v) ( x ) 
=
f ( v ( x ) )
:= a_3 v ^3(x)   + a_2 v ^2(x)  + a_1 v (x) + a_0,  
\
x \in (0, 1),
\
a_3 < 0,
\
a_2, a_1, a_0 \in \R,
\,
v \in L^{6} ( D; \R)
.
\end{equation*}
\end{assumption}
It is easy to find constants $L_0, L_1 \in (0, \infty)$ such that
\begin{equation}\label{eq:F-one-sided-condition}
\begin{split}
\langle u - v,  F (u) - F (v) \rangle & \leq L_0 \| u - v \|^2, \quad u , v \in V,
\\
\| F (u) - F (v) \| & \leq L_1 ( 1 + \| u \|_V^2 +  \| v \|_V^2 ) \| u - v \|,  \quad u, v \in V.
\end{split}
\end{equation}
The second property in \eqref{eq:F-one-sided-condition} immediately implies, for some $L_2 \in (0, \infty)$,
\begin{equation}
\| F (u) \|  \leq L_1 ( 1 + \| u \|_V^2 ) \| u  \| + \| F ( 0 ) \| \leq L_2 ( 1 + \| u \|_V^3 ),  \quad u \in V.
\end{equation}
Evidently, the constants $ L_0, L_1, L_2$ only depend on coefficients of $f$, i.e., $a_0, a_1, a_2, a_3$.
\begin{assumption}[Noise process]\label{ass:Phi}
Let $\{ W (t) \}_{t \in[0, T]}$ be a cylindrical $I$-Wiener process on a 
probability space  $\left(\Omega,\mathcal {F},\mathbb{P} \right)$ 
with a normal filtration $\{\mathcal{F}_t\}_{ t\in [0, T] }$, represented by a formal series,
\begin{equation} \label{eq:Wiener-representation}
W (t) := \sum_{n = 1}^{\infty} \beta_n ( t ) e_n,
\quad
t \in [0, T],
\end{equation}
where $\{ \beta_n ( t ) \}_{n \in \N}, t \in [0, T]$ is a sequence of independent real-valued standard Brownian motions 
and $\{e_n = \sqrt{2}\sin(n \pi x), \,x \in(0,1) \}_{n \in \N}$ is a complete orthonormal basis  of $H$.
%
%
\end{assumption}

\begin{assumption}[Initial value] \label{ass:X0}
Let the initial data $X_0 \colon  \Omega \rightarrow  H $, given by $ X_0 ( \cdot ) = u_0 ( \cdot ) $,
be an $\mathcal{F}_0/\mathcal{B}(H)$-measurable random variable.
For sufficiently large positive number $ p_0 \in \N $ and for any $\beta < \tfrac12$, 
{\color{black}{
there exists constants $K_\beta$ depending on $\beta$, $p_0$ 
and $\hat{K}_V$ only depending on $p_0$ such that}}
\begin{equation}
{\color{black}{
\| X_0 \|_{{L^{p_0}(\Omega, \dot{H}^\beta )} }  \leq K_{\beta} < \infty,
\qquad
 \| X_0 \|_{{L^{p_0}(\Omega, V )} }  \leq \hat{K}_V < \infty.
 }}
\end{equation}
\end{assumption}
%
%

%
At the moment, we are prepared to formulate the concrete problem \eqref{eq:SGL-concrete} as 
an abstract stochastic evolution equation in the Hilbert space $H$, 
\begin{equation}\label{eq:SGL-abstract}
\begin{split}
\left\{
    \begin{array}{lll} \dd X(t) + A X(t)\, \dd t = F ( X(t) ) \,\dd t +  \dd W(t), \quad  t \in (0, T], \\
     X(0) = X_0,
    \end{array}\right.
\end{split}
\end{equation}
where $X (t, \cdot ) = u (t, \cdot) $ and the abstract items $A, F, X_0$ 
are defined in Assumptions \ref{ass:A}-\ref{ass:X0}.
The above assumptions suffice to establish the well-posedness and regularity results of the mild solution 
to  \eqref{eq:SGL-abstract}, defined by \eqref{eq:mild-solution} later. Before that, we provide some 
regularity properties of the stochastic convolution.
\begin{lem}
Let Assumptions \ref{ass:A}, \ref{ass:Phi} be fulfilled.
For any $ p \in [ 2, \infty ) $
there exists {\color{black}{a constant $C_0$ dependent of $p$ but independent of $ \beta \in [0, \tfrac12) $}}  such that 
the stochastic convolution $ \{  \mathcal{O}_t \}_{t \in [0, T]}$ satisfies
\begin{align}
\| \mathcal{O}_t  \|_{L^p(\Omega, \dot{H}^\beta )} 
& \leq
C_0
{\color{black}{
\|A^{\frac{\beta-1}{2}}  \|_{\mathcal{L}_2 (H) }
}} 
< \infty,
\quad
\mbox{with}
\quad 
\mathcal{O}_t : = \int_0^t E(t-s) \dd W(s),
\\
%
%
%
\| \mathcal{O}_t - \mathcal{O}_s \|_{L^p(\Omega, H)} 
& \leq
C_0 
{\color{black}{
\|A^{\frac{\beta-1}{2}}  \|_{\mathcal{L}_2 (H) }
}}
(t -s)^{ \frac{\beta}{2} } , 
\quad 0 \leq s < t \leq T.
\label{eq:lem-stoch-conv-regularity}
\end{align}
Moreover, we have
\begin{equation}
\label{eq:lem-stoch-conv-spatial-regularity}
\E \Big[ \sup_{t \in [0, T] } \| \mathcal{O}_t \|_V^p   \Big]   < \infty.
\end{equation}
\end{lem}
Recalling \eqref{eq:semigroup-integral} one can validate the first two estimates easily, see, e.g., \cite[Theorem 2.4.]{wang2015note}.
The last assertion can, e.g., be found in \cite[Proposition 4.3]{da2004kolmogorov} and \cite[Lemma 6.1.2]{Cerrai2001second}.
Owing to the above regularity properties of the stochastic convolution, we can get the corresponding regularity properties of the
mild solution to \eqref{eq:SGL-abstract} as follows.
\begin{thm} \label{thm:SPDE-regularity-result}
Under Assumptions \ref{ass:A}-\ref{ass:X0}, SPDE \eqref{eq:SGL-abstract} possesses a unique mild solution 
$X: [0,T] \times \Omega \rightarrow V$ with continuous sample paths,
determined by, 
\begin{equation}\label{eq:mild-solution}
    X(t)
    =
    E(t) X_0
    +
    \int_0^t E(t-s) F ( X( s ) ) \, \dd s
    +
    \mathcal{O}_t
\quad  \mathbb{P} \mbox{-a.s.}
.
\end{equation}
For any $p \in [2, \infty)$ and $ \beta < \tfrac12$, {\color{black}{there exists positive constants $C_1, C_2$ 
depending on $p, T, \{a_i\}_{i=0}^{3}$  but not  depending on $\beta$}} such that,
\begin{align}
\label{eq:optimal.regularity1}
\sup_{ t \in [0, T] } \| X ( t ) \|_{L^p(\Omega, V )} 
        & \leq 
        C_1 \big( 1 + \| X_0 \|_{L^p(\Omega, V )}  \big),
\\
\label{eq:optimal.regularity2}
\sup_{ t \in [0, T] } \| X ( t ) \|_{L^p(\Omega, \dot{H}^{\beta } )} 
       & \leq 
       C_2 \big( 1 + \| X_0 \|_{L^ { p }(\Omega, \dot{H}^{\beta } ) } + \| X_0 \|_{L^ { 3p }(\Omega, V ) }^3 
         +
         {\color{black}{ \|A^{\frac{\beta-1}{2}}  \|_{\mathcal{L}_2 ( H ) }  }}
         \big).
\end{align}
Moreover, there exists a constant $C_3 \in [0, \infty)$ 
{\color{black}{depending on $p, T, \{a_i\}_{i=0}^{3}$ but not  depending on $\beta$}} such that,
for any $\beta < \tfrac12$ and $0 \leq s < t \leq T$,
\begin{equation}
\label{eq:optimal.regularity3}
\| X ( t ) - X ( s ) \|_{L^p(\Omega, H)} \leq C_3 
{\color{black}{
\big( 1 + \| X_0 \|_{L^ { p }(\Omega, \dot{H}^{\beta } ) } + \| X_0 \|_{L^ { 3p }(\Omega, V ) }^3 
         +
         \|A^{\frac{\beta-1}{2}}  \|_{\mathcal{L}_2 ( H ) }
         \big)
}}
(t -s)^{ \frac{\beta}{2} }.
\end{equation}
\end{thm}
The existence of the unique mild solution and the regularity assertion \eqref{eq:optimal.regularity1} are based on 
\cite[Proposition 6.2.2]{Cerrai2001second} and \eqref{eq:lem-stoch-conv-spatial-regularity}.
The rest of estimates in Theorem \ref{thm:SPDE-regularity-result} can be verified by standard arguments 
(consult, e.g., \cite[Theorem 2.4.]{wang2015note}).


\section{Spatial semi-discretization} \label{sect:spatial-discret}

This section concerns the error analysis for a spectral Galerkin spatial semi-discretization 
of the underlying problem \eqref{eq:SGL-abstract}.
For $N\in \mathbb{N}$ we define a finite dimensional subspace of $H$ by
\begin{equation}\label{eq:space.HN}
 H^N := \mbox{span} \{e_1, e_2, \cdots, e_N \},
\end{equation}
and the projection operator $P_N \colon \dot{H}^{\alpha}\rightarrow H^N$ by
$
  P_N \xi = \sum_{i=1}^N \langle \xi, e_i \rangle e_i, 
  \forall\, \xi \in \dot{H}^{\alpha}, \, \alpha \in \R.
$
Here $H^N$ is chosen as the linear space spanned by the $N$ first eigenvectors of 
the dominant linear operator $A$. It is not difficult to deduce that
\begin{equation}\label{eq:P-N-estimate}
\| ( P_N - I ) \varphi \|  \leq \lambda_{N+1}^{- \frac {\alpha}{2} } \|\varphi\|_{\alpha}
\leq N^{ - \alpha } \|\varphi\|_{\alpha}, 
\quad \forall \: \varphi \in \dot{H}^{\alpha}, \: \alpha \geq0.
\end{equation}
Additionally, define $A_N \colon H \rightarrow H^N$ as
$A_N = A P_N $, which generates an analytic semigroup $E_N(t) = e^{-t A_N}$, $t \in [0, \infty)$ in $H^N$. 
Then the spectral Galerkin approximation of \eqref{eq:SGL-abstract} results in the following finite dimensional SDEs,
\begin{equation}\label{eq:spectral-spde}
\begin{split}
\left\{
    \begin{array}{lll} \dd X^N(t) + A_N X^N ( t ) \, \dd t =  F_N ( X^N (t) ) \, \dd t  +  P_N  \, \dd W (t), 
    \quad t \in (0, T], \\
     X^N(0) = P_N X_0,
    \end{array}\right.
\end{split}
\end{equation}
where we write $ F_N : = P_N F $ for short. It is clear to see that \eqref{eq:spectral-spde} admits a unique solution in $H^N$.
By the variation of constant, the corresponding solution can be written as
\begin{equation}\label{eq:Spectral-Galerkin-mild}
    X^N(t)
    =
    E_N(t) P_N X_0
    +
    \int_0^t E_N (t - s ) P_N F ( X^N( s ) ) \, \dd s
    +
    \int_0^t E_N(t-s) P_N \, \dd W (s), \:  \mathbb{P} \mbox{-a.s.}.
\end{equation}
%
{\color{black}
{In the error analysis for the spatial semi-discretization \eqref{eq:spectral-spde}, 
we require $\| P_N X ( t ) \|_{ L^p(\Omega, V ) }, t > 0 $ to be
uniformly bounded with respect to $N \in \N$, which can be deduced based on the following two auxiliary results. 
The first one is a direct consequence of \cite[Lemma 5.4]{blomker2013galerkin} with $t_1 = 0$ and 
we follow the same arguments there to provide a short proof here.
}} 
\begin{lem} \label{Lem:PN-Stoch-Conv-V-bound}
Let $ \{  \mathcal{O}_t \}_{t \in [0, T]}$ be the stochastic convolution defined by \eqref{eq:lem-stoch-conv-spatial-regularity}
and let Assumptions \ref{ass:A}, \ref{ass:Phi} be fulfilled. Then
for any $p \in [2, \infty )$ it holds
\begin{equation}\label{eq:lemma-PN-Stoch-Conv-V-bound}
\sup_{ t \in [0, T], N \in \N }
\| \mathcal{O}_t^N \|_{L^p(\Omega, V ) } 
< \infty,
\qquad
\mbox{ with }
\,\,
\mathcal{O}_t^N := P_N \mathcal{O}_t.
\end{equation}
\end{lem}
{\color{black}{
{ \it Proof of Lemma \ref{Lem:PN-Stoch-Conv-V-bound}}.
Owing to Assumptions \ref{ass:A}, \ref{ass:Phi}, we perform the expansion of the stochastic convolution $O_t^{N}$ 
and use the It\^o isometry to get, for any $N \in \N, t \in [0, T], x, y \in D$,
\begin{align}\label{formula1}
\begin{split}
\E\big[\big| \mathcal O_t^N ( x )
   -
   \mathcal O_t^N ( y )\big|^2\big]
   &\leq
   \sum_{1 \leq i \leq N}
   \E\Big[
   \Big|
   \int_0^t e^{-\lambda_i(t- s )}
   \dd
   \beta_{i}(s)
   \Big|^2
   \Big]
   \big|e_i(x)-e_i(y)\big|^2
   \\&
   \leq
   \sum_{i\in \N}
   ( 2 \lambda_i )^{-1}
   ( \sqrt{2} \pi i )^{ \frac45 } | x - y |^{ \frac45}
   \big(
   |e_i(x)|+|e_i(y)|
   \big)^{ \frac65 }
   \\&
   \leq
   2^{\frac65} \pi^{-\frac65}
   | x - y|^{ \frac45 }
   \sum_{i\in \N}
   i^{-\frac65}
   \\&
   \leq 3
   | x - y|^{ \frac45 }.
\end{split}
\end{align}
In the same manner, one can acquire, for any $N \in \N, t \in [0, T]$,
\begin{align}
\sup_{x \in D}
\E
\big[
| \mathcal O_t^N ( x ) |^2
\big]
\leq 2.
\end{align}
Using the Sobolev embedding inequality $W^{\tfrac 15 , p} \subset V$, $p > 5$ and 
noting that the stochastic convolution $O_t^{N}$ is Gaussian one can deduce, for any $N \in \N, t \in [0, T]$,
\begin{equation}
\begin{split}
  \E \big[ \| \mathcal O_t^{N} \|_{V}^p \big]
 & \leq
   C  \int_0^1 \E \big[ | \mathcal O_t^{N}(x) |^p \big] \dd x
   +
   C \int_0^1 \int_0^1
   \frac{\E\big[\big| \mathcal O_t^{N}(x)
   -
   \mathcal O_t^{N}(y)\big|^p\big]}
   {|x-y|^{\tfrac {p}{5} + 1}}
   \, \dd x \dd y
\\
 &\leq
 C  \int_0^1
 \big(
 \E \big[ | \mathcal O_t^{N}(x) |^2 \big]
 \big)^{\tfrac p2} \dd x
   +
   C \int_0^1 \int_0^1
   \frac{\big(\E\big[\big| \mathcal O_t^{N}(x)
   -
   \mathcal O_t^{N}(y)\big|^2\big]\big)^{\tfrac p2}}
   {|x-y|^{\tfrac {p}{5} +1}}
   \, \dd x \dd y
  \\
  &
   \leq
   C \Big(
   1
   +
   \int_0^1 \int_0^1 \big| x-y \big|^{\tfrac {p}{5} -1}
   \dd x \dd y
   \Big)<\infty.
\end{split}
\end{equation}
Using H\"{o}lder's inequality yields the desired assertion in the case $p \in [2, 5]$. $\square$

The second one concerns the smoothing property of the analytic semigroup.
}
}
%
\begin{lem} \label{lem:semigroup-regularity0}
Let Assumptions \ref{ass:A} be fulfilled.
For any $ N \in \N $ and $ \psi \in \dot{H}^{\gamma},  \gamma \in [0, \tfrac12 )$, it holds that
\begin{equation} \label{eq:lem-V-H-regularity}
\| P_N E ( t ) \psi \|_V  \leq 
2^{ \gamma } 
\big ( 
\tfrac{ 5 - 4 \gamma }{ 2 \pi ( 1 - 2 \gamma ) } 
\big )^{ \frac12 }
t^{ \frac { 2 \gamma -1} {4} } \| \psi \|_{ \gamma },
\quad
t > 0,
\,
\gamma \in [0, \tfrac12 ).
\end{equation}
\end{lem}
{ \it Proof of Lemma \ref{lem:semigroup-regularity0}. }
Elementary facts readily yield
\begin{equation} 
\begin{split}
\| P_N E ( t ) \psi \|_V
& =
\sup_{ x \in [0, 1] }
\Big |
\sum_{ i = 1 }^N e^{ - \lambda_i t  } \langle \psi, e_i  \rangle e_i {\color{black}{(x)}}
\Big |
\leq
\sqrt{2}
\sum_{ i = 1 }^N  e^{ - \lambda_i t  } 
 | 
\langle \psi, e_i  \rangle
 |
\\
& \leq
\sqrt{2}
\Big(
\sum_{ i = 1 }^N  \lambda_i^{ - \gamma } e^{ - 2\lambda_i t  } 
\Big)^{1/2}
\Big(
\sum_{ i = 1 }^N
\lambda_i^{ \gamma }
 | 
\langle \psi, e_i  \rangle
 |^2
 \Big)^{1/2}
\\ &
\leq
 \sqrt{2} \pi^{ - \gamma }
\Big(
 \int_0^{ \infty }  x ^{ - 2 \gamma } e^{ - 2 \pi^2 x^2 t  }  \dd x
\Big)^{1/2}
\|
\psi
\|_{ \gamma }
\\ &
=
2^{ \gamma } \pi^{ - \frac12 } t^{ \frac{ 2 \gamma -1 }{ 4 } } 
\Big(
 \int_0^{ \infty }  y ^{ - 2 \gamma } e^{ - y^2/2  }  \dd y
\Big)^{1/2}
\|
\psi
\|_{ \gamma }
\\ &
\leq
2^{ \gamma } \big ( \tfrac{ 5 - 4 \gamma }{ 2 \pi ( 1 - 2 \gamma ) } \big )^{ \frac12 }
t^{ \frac { 2 \gamma -1} {4} } \| \psi \|_{ \gamma },
\end{split}
\end{equation}
as required.
$\square$

{\color{black}{
Equipped with the above two estimates, one can prove that $\| P_N X ( t ) \|_{ L^p(\Omega, V ) }, t > 0 $ is
uniformly bounded with respect to $N \in \N$ as follows.
}}
\begin{lem}  \label{lem:PN-Xt-V-bound}
Let $ \{ X ( t ) \}_{ t \in [0, T] } $ be the mild solution to \eqref{eq:SGL-abstract}, defined by \eqref{eq:mild-solution}.
Then for any $p \in [2, \infty )$ and $ \gamma \in [0, \tfrac12) $
{\color{black}{there exists a positive constant $C$,
depending on $T, p, \{a_i\}_{i=0}^{3}$, but not depending on $N, \beta, X_0$}} such that
\begin{equation}
\sup_{N \in \N } 
\| P_N X ( t ) \|_{ L^p(\Omega, V ) } 
\leq
{\color{black}{
2^{ \gamma } 
\big ( 
\tfrac{ 5 - 4 \gamma }{ 2 \pi ( 1 - 2 \gamma ) } 
\big )^{ \frac12 } 
\| X_0 \|_{ L^p(\Omega, \dot{H}^{ \gamma } ) } 
t^{ \frac {2 \gamma - 1} {4}}
+
C 
\big( 1 + \| X_0 \|^3_{L^{3p}(\Omega, V )}  \big)
}}
,
\quad
t \in ( 0, T].
\end{equation}
\end{lem}
{\it Proof of Lemma \ref{lem:PN-Xt-V-bound}. }
%
Observing that
$
E_N(t) P_N = E (t) P_N
$
and using Lemmas \ref{Lem:PN-Stoch-Conv-V-bound}, \ref{lem:semigroup-regularity0} show
\begin{align}
& \| P_N X(t) \|_{ L^p(\Omega, V ) }
    \leq
    \| E(t) P_N X_0 \|_{ L^p(\Omega, V ) }
    +
    \int_0^t \| E(t-s) F_N ( X( s ) ) \|_{ L^p(\Omega, V ) } \, \dd s
    +
    \| P_N \mathcal{O}_t \|_{ L^p(\Omega, V ) }
\nonumber
    \\
    & \quad \leq
    2^{ \gamma } 
\big ( 
\tfrac{ 5 - 4 \gamma }{ 2 \pi ( 1 - 2 \gamma ) } 
\big )^{ \frac12 }
    t^{\frac{2 \gamma - 1}{4}}  \| X_0 \|_{ L^p(\Omega, \dot{H}^{ \gamma } ) } 
    +
    \int_0^t  ( t - s )^{ - \frac{ 1 } { 4 } } \|  F ( X( s ) ) \|_{ L^p(\Omega, H ) } \, \dd s
    +
    \| P_N \mathcal{O}_t \|_{ L^p(\Omega, V ) }
\nonumber
    \\
    & \quad \leq
    2^{ \gamma } 
\big ( 
\tfrac{ 5 - 4 \gamma }{ 2 \pi ( 1 - 2 \gamma ) } 
\big )^{ \frac12 }
    t^{\frac{2 \gamma - 1}{4}}  \| X_0 \|_{ L^p(\Omega, \dot{H}^{ \gamma } ) } 
    +
    \tfrac{4}{3} t^{\frac34} \sup_{s \in [0, T]} \| F ( X( s ) ) \|_{ L^p(\Omega, H ) } 
    +
    \| P_N \mathcal{O}_t \|_{ L^p(\Omega, V ) }
\nonumber
    \\
    & \quad \leq
    2^{ \gamma } 
\big ( 
\tfrac{ 5 - 4 \gamma }{ 2 \pi ( 1 - 2 \gamma ) } 
\big )^{ \frac12 }
    t^{\frac{2 \gamma - 1}{4}}  
    \| X_0 \|_{ L^p(\Omega, \dot{H}^{ \gamma } ) } 
    +
    C \big ( 1 + \sup_{s \in [0, T]} \| X( s ) \|_{ L^{3p}(\Omega, V ) }^3 \big )
    +
    \| P_N \mathcal{O}_t \|_{ L^{p} (\Omega, V ) }.
\end{align}
Owing to \eqref{eq:optimal.regularity1}, \eqref{eq:lemma-PN-Stoch-Conv-V-bound} and Assumption \ref{ass:X0}, 
one can arrive at the expected estimate.
$\square$

Throughout this paper, 
by $C$ we mean deterministic constants, not necessarily the same at each occurrence, 
but independent of the discretization parameters.
Now we are prepared to carry out convergence analysis for 
the spectral Galerkin discretization \eqref{eq:spectral-spde}.
\begin{thm}[Spatial error estimate]
\label{thm:space-main-conv}
Let Assumptions \ref{ass:A}-\ref{ass:X0} hold. 
Let  $X(t)$ and $X^N(t)$ be defined through \eqref{eq:SGL-abstract} 
and \eqref{eq:Spectral-Galerkin-mild}, respectively. Then
for any $ \beta < \tfrac12 $, $ p \in [ 2, \infty )$ and $N \in \N$, 
{\color{black}{there exists a positive constant $C$,
depending on $T, p, \{a_i\}_{i=0}^{3}$, but not depending on $N, \beta, X_0$}} such that
\begin{equation}\label{eq:spatial-error}
\sup_{ t \in [ 0, T ] }
\| X(t) - X^N(t) \|_{L^p ( \Omega; H) } \leq
{\color{black}{
C 
\Big(
1 
+
 \| X_0 \|^{12}_{L^{12p}(\Omega, V )}
 +
 \| X_0 \|_{L^ { 4 p }(\Omega, \dot{H}^{\beta \vee [ (p-1)/(2p) ] } ) }^4
 +
 \|A^{\frac{\beta-1}{2}}  \|_{\mathcal{L}_2 }^2
\Big)
}}
N^{ - \beta }.
\end{equation}
\end{thm}
%
{\color{black}{Here and below we denote $a \vee b : = \max{ (a, b )}$.}} 
Clearly, the above convergence rate $\beta < \tfrac12$  can be arbitrarily close to $\tfrac12$ 
but  can not reach $\tfrac12$, since the quantity $ \|A^{\frac{\beta-1}{2}}  \|_{\mathcal{L}_2 }$,
as calculated in \eqref{eq:AQ_condition},  explodes when $\beta$ tends to $\tfrac12$. 
This comment also applies to the full approximation error estimates in section \ref{sect:full-discretization}.

{\it Proof of Theorem \ref{thm:space-main-conv}. } The triangle inequality along with \eqref{eq:optimal.regularity2},
\eqref{eq:P-N-estimate}  provides us that
\begin{equation}\label{eq:space-err-proof-split}
\begin{split}
\| 
 X(t) - X^N(t)
 \|_{L^p( \Omega, H)}
& \leq
\| 
 ( I - P_N ) X(t)
 \|_{L^p( \Omega, H)}
 +
 \| 
 P_N X(t) - X^N(t)
 \|_{L^p( \Omega, H)}
\\
& \leq
N^{ - \beta }
\| X( t ) \|_{L^{ p } ( \Omega; \dot{H}^{ \beta }) } 
+
 \| 
 e^N_t
 \|_{L^p( \Omega, H)}
\\ 
& \leq
{\color{black}{
C_2 \big( 1 + \| X_0 \|_{L^ { p }(\Omega, \dot{H}^{\beta } ) } + \| X_0 \|_{L^ { 3p }(\Omega, V ) }^3 
         +
         \|A^{\frac{\beta-1}{2}}  \|_{\mathcal{L}_2 } 
         \big) 
}}
N^{ - \beta }
\\
&
\quad
+
 \| 
 e^N_t
 \|_{L^p( \Omega, H)}
,
\end{split}
\end{equation}
where $ e^N_t : = P_N X(t) - X^N(t) = \int_0^t E_N ( t - s )  \big [ F_N ( X ( s ) ) - F_N ( X^N ( s ) ) \big] \, \dd s $ satisfies
\begin{equation}
\frac{\mbox{d} } { \mbox{d} t} e^N_t
=
- A_N e^N_t  +  F_N ( X(t) ) - F_N (X^N (t) )
=
- A e^N_t  +  F_N ( X(t) ) - F_N (X^N (t) ).
\end{equation}
Therefore, using \eqref{eq:F-one-sided-condition} and the Young inequality gives
\begin{equation} \label{eq:spatial-error-before-final}
\begin{split}
& \frac{\mbox{d} } { \mbox{d} t} \| e^N_t \|^p
=
p \| e^N_t \|^{p-2} \big \langle  e^N_t, 
- A e^N_t  +  F ( X(t) ) - F (X^N (t) )
\big\rangle
\\
&
\quad \leq
p \| e^N_t \|^{p-2} \big \langle  e^N_t, 
  F ( P_N X (t) ) - F (X^N (t) ) \big \rangle
+
p \| e^N_t \|^{p-2} \big \langle e^N_t,   F ( X (t) ) - F (P_N X(t) ) \big \rangle
\\
& 
\quad  \leq
L_0  p \| e^N_t \|^p
+ 
p \| e^N_t \|^{p-1} \big \| F ( X (t) ) - F (P_N X(t) ) \big \|
\\
&
\quad \leq
( L_0 p + p -1 )  \| e^N_t \|^p + \big \| F ( X (t) ) - F (P_N X(t) ) \big \|^p
\\
&
\quad \leq
( L_0 p + p -1 )  \| e^N_t \|^p 
+  
C \big (  1 + \| X (t) \|_V^{2 p} +  
\| P_N X (t) \|_V^{2 p} \big ) \| (I - P_N) X(t) \|^p.
\end{split}
\end{equation}
Choosing $ \gamma = \tfrac{ p - 1 } { 2 p } \in 
( \tfrac{ p - 2 } { 2 p },  \tfrac12)$ in Lemma \ref{lem:PN-Xt-V-bound} 
and also considering \eqref{eq:optimal.regularity1}, \eqref{eq:optimal.regularity2} and \eqref{eq:P-N-estimate} assure
\begin{align}
\E [ \| e^N_t \|^p ]
& \leq
C \! \int_0^ t \E [ \| e^N_s \|^p ]  
+
\E \big[
     \big ( 1 +  \| X (s) \|_V^{2p} +  
\| P_N X (s) \|_V^{2p} \big ) \| (I - P_N) X(s) \|^p
     \big] \, \dd s
\nonumber 
\\
&
\leq
C \! \int_0^ t \E [ \| e^N_s \|^p ]  
+
  \big(   
     1 +  \| X (s) \|_{L^{4p}(\Omega, V) }^{2p} +  
\| P_N X (s) \|_{L^{4p}(\Omega, V) }^{2p} \big ) 
\big\| (I - P_N) X(s) \big\|_{L^{2p}(\Omega, H) }^{p}
   \, \dd s
\nonumber
\\
&
{\color{black}{
\leq
C \! \int_0^ t \E [ \| e^N_s \|^p ]  \, \dd s
+
C 
N^{- p \beta }
\sup_{s \in [0, T]}
\| X (s) \|_{L^{2p}( \Omega, \dot{H}^\beta) }^p
}}
\nonumber
\\
&
\quad
{\color{black}{
\times
\int_0^t
\big (
1
+
 \| X_0 \|^{6p}_{L^{12p}(\Omega, V )}
+
s^{-\frac12} 
 \| X_0 \|^{2p}_{L^{4p}(\Omega, \dot{H}^{(p-1)/(2p)} )}
\big)
\dd s
}}
\nonumber
\\
&
{\color{black}{
\leq
C \! \int_0^ t \E [ \| e^N_s \|^p ]  \, \dd s
+
C
N^{- p \beta }
 \big(
 1 + \| X_0 \|_{L^ { p }(\Omega, \dot{H}^{\beta } ) } + \| X_0 \|_{L^ { 3p }(\Omega, V ) }^3 
         +
         \|A^{\frac{\beta-1}{2}}  \|_{\mathcal{L}_2 }
 \big)^p
 }}
 \nonumber
\\
&
\quad
{\color{black}{
\times
\big( 1 
+
 \| X_0 \|^{6p}_{L^{12p}(\Omega, V )}
 +
 \| X_0 \|^{2p}_{L^{4p}(\Omega, \dot{H}^{(p-1)/(2p)} )}
\big )
}}
 \nonumber
\\
&
{\color{black}{
\leq
C \! \int_0^ t \E [ \| e^N_s \|^p ]  \, \dd s
+
C
N^{- p \beta }
\big( 1 
+
 \| X_0 \|^{6}_{L^{12p}(\Omega, V )}
 +
 \| X_0 \|_{L^ { 4 p }(\Omega, \dot{H}^{\beta \vee [ (p-1)/(2p) ] } ) }^2
 +
 \|A^{\frac{\beta-1}{2}}  \|_{\mathcal{L}_2 }
\big )^{2p}.
}}
\label{eq:spatial-error-estimate-final}
\end{align}
After the use of the Gronwall inequality, one can derive from \eqref{eq:space-err-proof-split} the desired error bound.
$\square$

\section{Spatio-temporal full discretization}
\label{sect:full-discretization}

This section is devoted to error analysis of a spatio-temporal full discretization, 
done by a time discretization of the spatially discretized problem \eqref{eq:spectral-spde}.
For $M \in \N$ we construct a uniform mesh on $[0, T]$ with $\tau = \tfrac{T}{M}$ being the time stepsize,
and propose a spatio-temporal full discretization as,
\begin{equation}
\begin{split}
\label{eq:full.Tamed-AEE}
Y^{M,N}_{t_{m+1} } =& E_N(\tau) Y_{t_m}^{M,N} + \frac{ A_N^{-1} \big ( I - E_N( \tau ) \big ) F_N (Y_{t_m}^{M,N})  }
{  1 + \tau \| F_N(Y_{t_m}^{M,N}) \|   } 
+ \!\int_{t_m}^{t_{m+1}}\! E_N(t_{m+1}-s) P_N \mbox{d} W(s)
\end{split}  
\end{equation}
for $m = 0, 1, ..., M - 1$ and $ Y^{M,N}_0 = P_N X_0 $.
Equivalently, the full discretization \eqref{eq:full.Tamed-AEE} can be written by
$ Y^{M,N}_0 = P_N X_0 $ and for $m = 0, 1, ..., M - 1$,
\begin{equation}
\begin{split}
\label{eq:full.Tamed-AEE-Semigroup}
Y^{M,N}_{t_{m+1} } =& E_N(\tau) Y_{t_m}^{M,N} + \int_{t_m}^{t_{m+1}} \frac{ E_N(t_{m+1}-s)\, F_N (Y_{t_m}^{M,N})  }
{  1 + \tau \| F_N(Y_{t_m}^{M,N}) \|   } \,\mbox{d}s
+ \!\int_{t_m}^{t_{m+1}}\! E_N(t_{m+1}-s) P_N \mbox{d} W(s).
\end{split}  
\end{equation}
%
%
Here we invoke a taming technique in \cite{hutzenthaler2012strong,wang2013tamed,Tretyakov2013fundamental,hutzenthaler15MEMAMS} 
for ordinary SDEs, and construct a nonlinearity-tamed accelerated exponential Euler (AEE) scheme as \eqref{eq:full.Tamed-AEE}.
The so-called AEE scheme without taming is originally introduced in \cite{jentzen2009overcoming},
to strongly approximate nearly linear parabolic SPDEs. Since the stochastic convolution is 
Gaussian distributed and diagonalizable on $\{  e_i \}_{i \in \N}$,  the scheme is much easier to simulate than 
it appears at first sight (see comments in section \ref{sec:numerical-result} for the implementation).
When the nonlinearity grows super-linearly,  one can in general not expect that  
the usual AEE schemes \cite{jentzen2011efficient,wang2015note, 
wang2014higher,qi2017accelerated,lord2016modified,jentzen2009overcoming} converge strongly,
based on the observations from 
\cite{hutzenthaler2011strong,beccari2019strong}.  
Also, we mention that analyzing the strong convergence rate is much more difficult than that in the finite dimensional SDE setting. 

{\color{black}{When analyzing strong convergence rates, a crucial element is to derive uniform moment bounds for the 
spatio-temporal full discretization. To do this, a key ingredient lies on transforming 
the continuous time extension \eqref{eq:full-discrete-continous} of
the full discretization \eqref{eq:full.Tamed-AEE-Semigroup} into a random differential equation \eqref{eq:Numerical-RPDE},
with a random perturbation appearing in the cubic nonlinear term. So the forthcoming subsection
attempts to reveal that the $V$-norm of the solution to the perturbed differential equation can be controlled by a norm of the perturbation 
(see Lemma \ref{lem:vN-controlled-ZN}). Based on this finding, one can further use the structure of the taming scheme
and a certain bootstrap argument (Lemma \ref{Lem:Full-discrete-MB} and \eqref{eq:moment-bound-part2}) 
to provide a priori moment bound of the approximation, 
as shown by Theorem \ref{thm:Numer-Moment-Bound} in subsection \ref{subsec:a-priori-moment-numerical}.  
Reecovering the higher temporal convergence  rate of order almost $\tfrac12$ also heavily hinges on 
commutativity properties of the nonlinearity (Lemma \ref{Lem:F-negative-soblev}) and 
the improved temporal H\"{o}lder regularity of the solution in negative Sobolev space 
(Lemma \ref{Lem:regularity-negative-soblev}), as presented in
subsection \ref{subset:negative-sobleve}. Armed with these useful results, the expected convergence rate 
is finally proved in subsection \ref{subset:full-error-analysis}.
}}
%
%
{\color{black}{\subsection{Useful estimates for a perturbed differential equation}
}}
\label{subsec:estimates-perturbed-PDE}
%
%
%
%
%
{\color{black}{In the first part, we restrict ourselves to the following perturbed differential equation in $H^N$, $ N \in \N $,
\begin{equation} \label{eq:determ-perturbed-PDE}
\left\{
    \begin{array}{ll}
    \frac{\partial v^N }{\partial t } = - A_N  v^N + P_N F ( v^N + z^N ),
    \quad t \in (0, T],
     \\
     v^N ( 0 ) = 0,
    \end{array}\right.
\end{equation}
where $F$ comes from Assumption \ref{ass:F} and $ z^N , v^N \colon  [0, T] \rightarrow H^N $. 
In what follows, we aim to show that the $V$-norm of the solution to the perturbed problem \eqref{eq:determ-perturbed-PDE} 
can be controlled by the $\mathbb{L}^9$-norm of the perturbation $z^N$ in the nonlinear mapping (see Lemma \ref{lem:vN-controlled-ZN}), 
which plays a crucial role in deducing uniform moment bounds for the full discretization \eqref{eq:full.Tamed-AEE}. 
To show this, we need further smoothing properties of the analytic semigroup $E(t), t > 0$ described as follows.
}
}
\begin{lem}\label{lem:semigroup-norm}
Let $t > 0$ and let $P_N, E (t)$ be defined as in the above sections. Then it holds that
\begin{equation} \label{eq:lem-semigroup-regularity}
\begin{split}
{\color{black}{
\| P_N E ( t ) \psi \|_{L^p( D )} 
}}
& 
{\color{black}{
\leq t^{ - \frac{p-2}{4p} } \| \psi \|,
\quad
\forall p \geq 2,
}} 
\\
\|  P_N E ( t ) \psi \|_V & \leq ( \tfrac{ t } { 2 } ) ^{ - \frac12 } \| \psi \|_{ L^1( D ) },
\\
\| P_N E ( t ) \psi \|_{ L^3( D ) } & \leq ( \tfrac{ t } { 2 } )^{ - \frac{5}{24} } \| \psi \|_{ L^{\frac43}( D ) }.
\end{split}
\end{equation}
\end{lem}
{ \it Proof of Lemma \ref{lem:semigroup-norm}. }  
{\color{black}{
Recalling \eqref{eq:lem-V-H-regularity} with $\gamma = 0$ helps us to infer,
\begin{equation}
\| P_N E ( t ) \psi \|_{ L^p (D) }^p \leq \| P_N E ( t ) \psi \|^2 \cdot \| P_N E ( t ) \psi \|_{ V }^{p-2}
\leq 
\| \psi \|^2 \cdot t^{ - \frac{p-2}{4} } \| \psi \|^{p-2} 
= 
t^{ - \frac{p-2}{4} } \| \psi \|^p,
\end{equation}
}}
which validates the first assertion. To arrive at the second one, 
we again use \eqref{eq:lem-V-H-regularity} with $\gamma = 0$ to get
\begin{equation}
\begin{split}
\| P_N E ( t ) \psi \|_V & = \| P_N E ( \tfrac{t}{2} ) P_N E ( \tfrac{t}{2} )  \psi \|_V 
\leq
( \tfrac{t}{2} )^{ - \frac14 } \| P_N E ( \tfrac{t}{2} )  \psi \|
= 
( \tfrac{t}{2} )^{ - \frac14 } \sup_{ \| \phi \| \leq 1 } \big | \langle P_N E ( \tfrac{t}{2} )  \psi,  \phi  \rangle \big |
\\ &
=
( \tfrac{t}{2} )^{ - \frac14 } \sup_{ \| \phi \| \leq 1 } \big | \langle   \psi,  P_N E ( \tfrac{t}{2} ) \phi  \rangle \big |
\leq
( \tfrac{t}{2} )^{ - \frac14 } \sup_{ \| \phi \| \leq 1 } \| \psi \|_{ L^1(D) } \cdot \| P_N E ( \tfrac{t}{2} ) \phi \|_V
\\ &
\leq
( \tfrac{t}{2} )^{ - \frac12 } \| \psi \|_{ L^1(D) }.
\end{split}
\end{equation}
Concerning the last inequality, one can similarly acquire
\begin{equation}
\begin{split}
\| P_N E ( t ) \psi \|_{ L^3( D ) } & \leq ( \tfrac{t}{2} )^{ - \frac{1}{12} } \| P_N E ( \tfrac{t}{2} )  \psi \|
= 
( \tfrac{t}{2} )^{ - \frac{1}{12} } \sup_{ \| \phi \| \leq 1 } \big | \langle P_N E ( \tfrac{t}{2} )  \psi,  \phi  \rangle \big |
\\
& =
( \tfrac{t}{2} )^{ - \frac{1}{12} } \sup_{ \| \phi \| \leq 1 } \big | \langle \psi,  P_N E ( \tfrac{t}{2} ) \phi  \rangle \big |
\leq
( \tfrac{t}{2} )^{ - \frac{1}{12} } \sup_{ \| \phi \| \leq 1 } \| \psi \|_{ L^{\frac43} (D) } 
\| P_N E ( \tfrac{t}{2} ) \phi  \|_{ L^{4} ( D ) }
\\
& \leq
( \tfrac{t}{2} )^{ - \frac{5}{24} }  \| \psi \|_{ L^{\frac43} (D) }.
\end{split}
\end{equation}
The proof is now completed. $\square$
%

It is easy to see, the perturbed problem \eqref{eq:determ-perturbed-PDE} has a unique solution in $ H^N $, which can be expressed by
\begin{equation}
v^N ( t ) = \int_0^t E_N (t - s ) P_N F ( v^N (s) + z^N (s) ) \, \dd s.
\end{equation}
Define norms $\| u \|_{ \mathbb{ L }^q( D \times [0, t]  ) } : = 
\big ( \int_0^t \| u (s) \|_{L^{q} (D) }^q \, \dd s  \big)^{\frac1q}, \, q \geq 1, \, t \in [ 0, T ] $. For the particular case $q = 2$, 
$\mathbb{ L }^q( D \times [0, t]  ) $ ($\mathbb{ L }^q $ for brevity) becomes a  Hilbert space with 
$\langle u, v \rangle_{ \mathbb{ L }^2 ( D \times [0, t]  ) } : = \int_0^t \langle u (s), v (s) \rangle \dd s $.
{\color{black}{The next lemma asserts that the $V$-norm of the solution to \eqref{eq:determ-perturbed-PDE} 
can be controlled by the $\mathbb{L}^9$-norm of the perturbation $z^N$, which will be essentially used 
in proving moment bounds of the approximations.
}}
\begin{lem} \label{lem:vN-controlled-ZN}
Let $ v^{N}, N \in \N $ be the solution to \eqref{eq:determ-perturbed-PDE}. For any $ t \in [ 0, T ] $, 
there exists a positive constant $C$, {\color{black}{dependent of $T$ but independent of $N$}}, such that
\begin{equation}
\label{eq:vN-controlled-ZN}
\| v^{N} ( t ) \|_V 
\leq 
C ( 1 + \| z^N \|_{ \mathbb{L}^9 ( D \times [0, t] ) }^9 ),
\quad
\forall \ t \in [0, T].
\end{equation}
\end{lem}
{ \it Proof of Lemma \ref{lem:vN-controlled-ZN}. }
The assertion is trivial for $ t = 0$. So we always suppose $ t > 0$ in the following. 
We divide the proof into two steps. 

{\it Step 1.} For any fixed $t \in ( 0, T]$, we claim first that, 
by setting $ \varrho_t :=  5 t^{\frac14} \max \{ \tfrac{ | a_2 | } { | a_3 | } , \big|  \tfrac{ a_1  } { a_3 } \big|^\frac12, 
\big | \tfrac {  a_0 } {  a_3  } \big|^{ \frac13 } \} $,
\begin{equation} \label{eq:claim1}
\| v^N \|_{ \mathbb{ L }^4( D \times [0, t] ) } 
\leq 
5  \| z^N \|_{ \mathbb{ L }^4( D \times [0, t] ) } 
\quad
\text{ or }
\quad
\| v^N \|_{ \mathbb{ L }^4( D \times [0, t] ) }   \leq \varrho_t.
\end{equation}
By deterministic calculus and noting  $ A_N v^N = A v^N$ for any $v^N \in H^N$,  we derive
\begin{equation}\label{eq:inner-produc-positive}
\begin{split}
0 
& \leq 
\tfrac12 \| v^N (t) \|^2
= 
 \int_0^t \left \langle v^N ( s ),  
- A_N v^N ( s ) 
+ 
P_N F ( v^N(s) + z^N (s)  ) \right\rangle \dd s
\\
& \leq
 \int_0^t \left \langle v^N ( s ),   
F ( v^N(s) + z^N (s)  ) \right\rangle \dd s
=
\langle v^N, F ( v^N + z^N ) \rangle_{ \mathbb{ L }^2 }.
\end{split}
\end{equation}
Noticing that $ a_3 < 0 $, for any $ v , z \in \R$,
\begin{equation}
\begin{split}
v f ( v + z ) & = a_3 v ( v + z )^3 + a_2 v ( v + z )^2 + a_1 v ( v + z ) + a_0 v
\\ &
\leq
a_3 v^4 + 3 a_3 v^3 z + a_3 v z^3 + a_2 v^3 + 2 a_2 v^2 z + a_2 v z^2 + a_1 v^2 + a_1 v z + a_0 v.
\end{split}
\end{equation}
After using the fact $a_3 < 0$ and the H\"{o}lder inequality, one derives
\begin{equation}
\begin{split}
\langle v^N, F ( v^N + z^N ) \rangle_{ \mathbb{ L }^2 } 
& \leq
a_3 \| v^N \|_{ \mathbb{ L }^4 }^4 
+
3 | a_3 | \| v^N \|_{ \mathbb{ L }^4 }^3 \| z^N \|_{ \mathbb{ L }^4 }
+
| a_3 | \| v^N \|_{ \mathbb{ L }^4 } \| z^N \|_{ \mathbb{ L }^4 }^3
\\
& \quad + 
| a_2 | t^{ \frac14 }  \| v^N \|_{ \mathbb{ L }^4 }^3
+ 
2 | a_2 | t^{ \frac14 } \| v^N \|_{ \mathbb{ L }^4 }^2 \| z^N \|_{ \mathbb{ L }^4 }
+
| a_2 |  t^{ \frac14 } \| v^N \|_{ \mathbb{ L }^4 } \| z^N \|_{ \mathbb{ L }^4 }^2
\\
&
\quad
+ 
| a_1 |  t^{ \frac12 } \| v^N \|_{ \mathbb{ L }^4 }^2
+
| a_1 |  t^{ \frac12 } \| v^N \|_{ \mathbb{ L }^4 }  \| z^N \|_{ \mathbb{ L }^4 }
+
| a_0 | t^{ \frac34 }   \| v^N \|_{ \mathbb{ L }^4 } .
\end{split}
\end{equation}
Assume the claim \eqref{eq:claim1} is false, namely, 
$
\| z^N \|_{ \mathbb{ L }^4( D \times [0, t] ) } < \frac{ 1 }{ 5 } \| v^N \|_{ \mathbb{ L }^4( D \times [0, t] ) }  
$
and
$
\| v^N \|_{ \mathbb{ L }^4( D \times [0, t] ) }   >  \varrho_t.
$
This enables us to derive, in the case $ \varrho_t > 0 $, i.e., $| a_0 | + | a_1 | + | a_2 | > 0$,
\begin{equation} \label{eq:contradiction}
\begin{split}
\langle v^N, F ( v^N + z^N ) \rangle_{ \mathbb{ L }^2 } 
& <
\big ( a_3 + \tfrac{ 3 | a_3 | }{ 5 } + \tfrac{ | a_3 | }{ 125 }  \big) \| v^N \|_{ \mathbb{ L }^4  }^4
+
\big( | a_2 | t^{ \frac14 } +  \tfrac{ 2 | a_2 |  t^{ \frac14 } }{ 5 } +  \tfrac{ | a_2 |  t^{ \frac14 } }{ 25 }  \big) \| v^N \|_{ \mathbb{ L }^4  }^3
\\
& \quad + 
\big ( |a_1| t^\frac12 +  \tfrac{ | a_1 |  t^{ \frac12 } }{ 5 } \big) \| v^N \|_{ \mathbb{ L }^4  }^2
+
| a_0 | t^{ \frac34 }   \| v^N \|_{ \mathbb{ L }^4 } 
\\
& <
\Big( 
a_3 + \tfrac{ 76 | a_3 | }{ 125 } 
+
\tfrac{ | a_2 | t^{ \frac14 } }{ \varrho_t } +  \tfrac{ 2 | a_2 |  t^{ \frac14 } }{ 5 \varrho_t } +  \tfrac{ | a_2 |  t^{ \frac14 } }{ 25 \varrho_t } 
+
\tfrac{ |a_1| t^\frac12 } { \varrho_t^2 } +  \tfrac{ | a_1 |  t^{ \frac12 } }{ 5 \varrho_t^2 }
+
\tfrac{ | a_0 | t^{ \frac34 } } { \varrho_t^3 }
\Big) 
\| v^N \|_{ \mathbb{ L }^4  }^4
\\
& \leq
{\color{black}{
\big( 
a_3 + \tfrac{ 76 | a_3 | }{ 125 } 
+
\tfrac{ | a_3 | }{ 5 } +  \tfrac{ 2 | a_3 | }{ 25 } +  \tfrac{ | a_3 |   }{ 125 } 
+
\tfrac{ | a_3 | }{ 25 } +  \tfrac{ | a_3 | }{ 125 }
+
\tfrac{ | a_3 | }{ 125 }
\big) 
\| v^N \|_{ \mathbb{ L }^4  }^4
}}
\\
& =
\tfrac{ 6 a_3 }{ 125 } \| v^N \|_{ \mathbb{ L }^4  }^4
< 0,
\end{split}
\end{equation}
{\color{black}{where we set $ \varrho_t :=  5 t^{\frac14} \max \{ \tfrac{ | a_2 | } { | a_3 | } , \big|  \tfrac{ a_1  } { a_3 } \big|^\frac12, 
\big | \tfrac {  a_0 } {  a_3  } \big|^{ \frac13 } \} $ such that all coefficients of  $\| v^N \|_{ \mathbb{ L }^4  }^4$ in the second step 
only consist of $|a_3|$ or $a_3$ and the resulting collection turns out to be negative.}}
The conclusion \eqref{eq:contradiction} thus contradicts \eqref{eq:inner-produc-positive}. 
When $ a_0  = a_1  = a_2 = 0$, i.e., $\varrho_t = 0$,  one can similarly derive
$
\langle v^N, F ( v^N + z^N ) \rangle_{ \mathbb{ L }^2 } 
<
\big ( a_3 + \tfrac{ 3 | a_3 | }{ 5 } + \tfrac{ | a_3 | }{ 125 }  \big) \| v^N \|_{ \mathbb{ L }^4  }^4
=
\tfrac{ 49 a_3 }{ 125 } \| v^N \|_{ \mathbb{ L }^4  }^4
\leq 0,
$
also contradicting \eqref{eq:inner-produc-positive}.
Therefore, the claim \eqref{eq:claim1} must be true.

{\it Step 2.}
Apparently, \eqref{eq:claim1} implies
\begin{equation}
\| v^N \|_{ \mathbb{ L }^4 (  D \times [0, t]   ) } 
\leq 
5  \| z^N \|_{ \mathbb{ L }^4 ( D \times [0, t] ) } 
+ \varrho_T,
\quad
\forall t \in [0, T].
\end{equation}
This together with the last inequality in \eqref{eq:lem-semigroup-regularity}, the property of the cubic nonlinearity 
and the H\"{o}lder inequality yields, for any $ t \in [ 0, T]$,
\begin{equation}\label{vN-L3-estimate}
\begin{split}
\| v^N ( t ) \|_{L^3 ( D ) }  
& 
\leq  
\int_0^t \| E (t - s ) P_N F ( v^N (s) + z^N (s) ) \|_{L^3 ( D ) } \, \dd s
\\
& \leq
\int_0^t ( \tfrac { t - s} {2} )^{ - \frac{5}{24} } \|  F ( v^N (s) + z^N (s) ) \|_{L^{\frac43} ( D ) } \, \dd s
\\
& \leq
C 
\int_0^t ( \tfrac { t - s} {2} )^{ - \frac{ 5 } { 24 } } \big( 1 + \| v^N (s) \|_{ L^4 (D) }^3 + \| z^N (s) \|_{ L^4 (D) }^3 \big) \, \dd s
\\
& \leq
C 
\Big (  \int_0^t ( \tfrac { t - s} {2} )^{ - \frac{ 5 } { 6 } } \, \dd s \Big )^{ \frac14 } 
\Big (  \int_0^t   1 + \| v^N (s) \|_{ L^4 (D) }^4 + \| z^N (s) \|_{ L^4 (D) }^4  \, \dd s \Big )^{ \frac34 } 
\\
& \leq
C \big ( 1 + \| z^N \|_{ \mathbb{L}^4 ( D \times [0, t] ) }^3  \big).
\end{split}
\end{equation}
Likewise, by virtue of the second inequality in \eqref{eq:lem-semigroup-regularity} instead, one obtains, 
for any $ t \in [ 0, T]$,
\begin{equation}
\begin{split}
\| v^N ( t ) \|_V  
& \leq  
\int_0^t \| P_N E ( t-s) F ( v^N (s) + z^N (s) ) \|_V \, \dd s
\\
& \leq
C 
\int_0^t ( \tfrac{ t - s}{2} )^{-\frac12} \| F ( v^N (s) + z^N (s) ) \|_{L^1(D)} \, \dd s
\\
& \leq
C 
\int_0^t ( t - s )^{-\frac12}  ( 1 + \| v^N (s) \|_{L^3(D)}^3 + \| z^N (s) \|_{L^3(D)}^3 )  \, \dd s
\\
& \leq
C \bigg ( 1 + \| z^N \|_{ \mathbb{L}^4 ( D \times [0, t] ) }^9 +  \int_0^t  ( t - s )^{-\frac12} \| z^N (s) \|_{L^3(D)}^3 \, \dd s   \bigg)
\\
& \leq 
C ( 1 + \| z^N \|_{ \mathbb{L}^9 ( D \times [0, t] ) }^9 ).
\end{split}
\end{equation}
The proof of Lemma \ref{lem:vN-controlled-ZN} is thus finished. $\square$

{\color{black}{
Before closing this subsection, we would like to add some comments on the proof of Lemma \ref{lem:vN-controlled-ZN}.
The key step is to show either 
$
\| v^N \|_{ \mathbb{ L }^4( D \times [0, t] ) } 
\leq 
\kappa  \| z^N \|_{ \mathbb{ L }^4( D \times [0, t] ) }
$ 
or
$
\| v^N \|_{ \mathbb{ L }^4( D \times [0, t] ) }   \leq \varrho_t
$, 
with some positive constants $\kappa$, $\varrho_t$ properly chosen.
As shown above, its proof is based on the reduction to absurdity and 
we choose $\kappa = 5$ and $ \varrho_t :=  5 t^{\frac14} \max \{ \tfrac{ | a_2 | } { | a_3 | } , \big|  \tfrac{ a_1  } { a_3 } \big|^\frac12, 
\big | \tfrac {  a_0 } {  a_3  } \big|^{ \frac13 } \} $ in \eqref{eq:contradiction} to result in a contradiction. 
We mention that the choice of $\kappa$, $\varrho_t$ is, however, not unique.
}}

\subsection{A priori moment bounds of the approximations}
\label{subsec:a-priori-moment-numerical}
{\color{black}{
The aim of this subsection is to prove a priori moment bounds for the fully discrete approximation \eqref{eq:full.Tamed-AEE}
(or \eqref{eq:full.Tamed-AEE-Semigroup} equivalently), 
which is essentially based on the estimate \eqref{eq:vN-controlled-ZN} obtained in the previous subsection 
as well as a certain bootstrap argument.}}
First, we define 
%
\begin{equation}
\lfloor t \rfloor : =  t_i = i \tau,  \quad \text{ for } t \in [t_i, t_{i+1}), 
\quad 
i  \in \{ 0, 1, ... ,M-1\} ,
\end{equation}
and introduce a continuous version of the fully discrete scheme \eqref{eq:full.Tamed-AEE-Semigroup} as, 
\begin{equation} \label{eq:full-discrete-continous}
\begin{split}
Y^{M,N}_t 
= & 
E_N(t) Y_0^{M,N} 
+ 
\int_{0}^{t} \tfrac{ E_N(t - s)\, F_N(Y^{M,N}_{\lfloor s \rfloor})  }
{  1 + \tau \| F_N(Y^{M,N}_{\lfloor s \rfloor}) \|   } \,\mbox{d}s 
+ 
\mathcal{O}_t^N, 
\quad
\,
t \in [0, T].
\end{split}
\end{equation}
%
Here we recall that $\tau = \tfrac{T}{M}$ is the uniform time stepsize. 
By $B^c$ and $\1_B$, we denote the complement and indicator function of a set $B$, respectively.
Additionally, we introduce a sequence of decreasing subevents
\begin{equation}\label{Eq:tame.Omega}
  \Omega_{R,t_i}:= 
  \Big \{ 
  \omega \in \Omega: \sup_{j  \in \{ 0, 1, ... ,i\} } \| Y^{M,N}_{ t_j }(\omega) \|_V \leq R \Big\},
  \quad
  R \in (0, \infty), \, i  \in \{ 0, 1, ... ,M\}.
\end{equation}
It is clear that $ \1_{\Omega_{R,t_i}} \in \F_{t_i}$ for $ i  \in \{ 0, 1, ... ,M\} $ and 
$\1_{ \Omega_{R,t_i} } \leq \1_{ \Omega_{R, t_j} }$
for $t_i \geq t_j$ since $ \Omega_{R,t_i} \subset  \Omega_{R, t_j}, t_i \geq t_j $.
Furthermore, we put additional assumptions on the initial data.
\begin{assumption}\label{ass:X0-additional}
For sufficiently large positive number $ p_0 \in \N $, the initial data $X_0$ obeys
\begin{equation} \label{eq:ass-X0-additional}
\sup_{ N \in \N } \| P_N X_0 \|_{L^{p_0}( \Omega, V)} < \infty.
\end{equation}
\end{assumption}
Thanks to the Sobolev embedding inequality, \eqref{eq:ass-X0-additional} is fulfilled provided $\| P_N X_0 \|_{L^{p_0}( \Omega, \dot{H}^{\gamma})} < \infty$ for any $\gamma > \tfrac12$. In what follows we start the bootstrap argument, by showing
$ \E\big[ \1_{ \Omega_{ R^{ \tau },t_m} } 
\| Y^{M,N}_{t_m} \|_V^{p} \big] < \infty$ and $    \E
   \big [
   \1_{\Omega^c_{R^{\tau},t_m}} \|Y_{t_m}^{M,N} \|_V^p
   \big ] < \infty $
for  subevents $ \Omega_{ R^{ \tau },t} $  with $ R^{ \tau } $ depending on $\tau$ carefully chosen.
\begin{lem}\label{Lem:Full-discrete-MB}
Let $ p \in [2, \infty)$ and $ R^{\tau} := \tau^{-\frac{\beta}{ 4 }} $ for any $ \beta \in (0, \tfrac12) $. 
Under Assumptions \ref{ass:A}-\ref{ass:X0}, \ref{ass:X0-additional}, 
the approximation process $Y^{M,N}_{t_i}, i \in \{ 0, 1,..., M \}$  produced by \eqref{eq:full.Tamed-AEE} obeys
\begin{equation} \label{eq:lem-moment-bound0}
\sup_{ M, N \in \N} \sup_{ i \in \{ 0, 1,..., M \}} 
\E\big[ \1_{\Omega_{ R^{ \tau }, t_{i-1} } } \| Y^{M,N}_{t_i} \|_V^{p} \big]
<
\infty,
\end{equation}
where we set $ \1_{\Omega_{R^{\tau},t_{-1} }} = 1$.
\end{lem}
{\it Proof of Lemma \ref{Lem:Full-discrete-MB}. }
 The proof heavily relies on the use of Lemma \ref{lem:vN-controlled-ZN}.
In order to apply it, we introduce a process $Z_t^{ M, N} $ given by,
\begin{equation} \label{eq:Zt-defn}
\begin{split}
Z_t^{ M, N}
 & :=  
 E_N(t) Y_0^{M,N} 
 +
   \int_0^t E_N ( t - s )
   \Big[
   \tfrac{ P_N F( Y^{M,N}_{\lfloor s \rfloor} ) }{ 1 + \tau \| P_N F( Y^{M,N}_{\lfloor s \rfloor} ) \| }
   -
   P_N F( Y^{M,N}_s )
   \Big]
   \,\dd s
 +
 \mathcal{O}_t^N
 \\
 & = 
 E_N(t) Y_0^{M,N}
 +
   \int_0^t E(t-s) P_N
   \big [
   F(Y^{M,N}_{\lfloor s \rfloor})
   -
   F(Y^{M,N}_s)
   \big ]
   \,\dd s
\\
 &
 \quad +
   \int_0^t E(t-s)
   \Big [
   \tfrac{ P_N F ( Y^{M,N}_{ \lfloor s \rfloor } ) } {1 + \tau\| P_N F(Y^{M,N}_{\lfloor s \rfloor}) \|}
   -
   P_N F ( Y^{M,N}_{ \lfloor s \rfloor} )
   \Big ]
   \,\dd s
   +
   P_N \mathcal{O}_t,
 \quad
 t \in [0, T].
\end{split}
\end{equation}
Bearing this in mind, one can rewrite \eqref{eq:full-discrete-continous} as
\begin{equation} \label{eq:full-discrete-continous-form2}
\begin{split}
Y^{M,N}_t
& =
   \int_0^t E_N (t - s) P_N F ( Y^{M,N}_s ) \,\dd s
   +
   Z_t^{ M, N}.
\end{split}
\end{equation}
Further, we define $\bar{Y}^{M,N}_t$ as
\begin{equation} \label{eq:defn-Ybar}
   \bar{Y}^{M,N}_t
   :=
   Y^{M,N}_t - Z^{M,N}_t,
   \quad
   \text{ with }
   \quad
   \bar{Y}^{M,N}_0
   =
   0
   .
\end{equation}
So, we recast \eqref{eq:full-discrete-continous-form2} as
\begin{equation}
   \bar{Y}^{M,N}_t
   =
   \int_0^t E_N (t-s) P_N F ( \bar{Y}^{M,N}_s+Z^{M,N}_s )
   \,\dd s,
   \quad
   t \in [0, T],
\end{equation}
which satisfies
\begin{equation} \label{eq:Numerical-RPDE}
   \frac{\dd  } { \dd t} \bar{Y}^{M,N}_t
   =
   - A_N \bar{Y}^{M,N}_t
   +
   P_N F(\bar{Y}^{M,N}_t+Z^{M,N}_t),
   \quad
   t \in ( 0, T],
 \quad
   \bar{Y}^{M,N}_0
   =
   0.
\end{equation}
Now one can employ Lemma \ref{lem:vN-controlled-ZN} to obtain, 
\begin{equation}
\| \bar{Y}^{M,N}_t \|_V \leq C ( 1 + \| Z^{M,N} \|_{ \mathbb{L}^9 ( D \times [0, t] ) }^9 ),
\quad
t \in [0, T],
\end{equation}
where $Z^{M,N}_\cdot$ is defined by \eqref{eq:Zt-defn}.
Thus, for any $ i \in \{ 0, 1,...,M \} $,
\begin{equation}\label{eq:Zt-Yt}
\begin{split}
\E \big[ \1_{\Omega_{R,t_{i-1}}} \| \bar{ Y }^{ M, N }_{t_i} \|_V^p \big]
& \leq 
C 
\big( 1 + \E \big[ \1_{\Omega_{R,t_{i-1} }} \| Z^{M,N} \|_{ \mathbb{L}^{9p}( D \times [0, t_i] )  }^{9p}  \big] 
\big)
\\
& \leq 
C 
\Big( 1 + \E \Big[ 
\1_{\Omega_{R,t_{i-1} }}  \!
\int_0^{t_i} \| Z^{M,N}_s \|_V^{9p}  \, \dd s \Big] 
\Big),
\end{split}
\end{equation}
where, for $s \in [0, t_i]$, $ i \in \{ 0, 1,...,M \} $, it holds that
\begin{equation} \label{eq:Z-V-subevents}
\begin{split}
  \1_{\Omega_{R,t_{i-1} }}
   \| Z^{M,N}_s \|_V
&
   \leq
   \| E_N( s ) Y_0^{M,N} \|_V
   +
   \1_{\Omega_{R,t_{i-1} }}
   \bigg\|
   \!
   \int_0^s 
   E ( s - r ) P_N
   \big [
   F( Y^{M,N}_r )
   -
   F( Y^{ M,N }_{ \lfloor r \rfloor } )
   \big ]
   \,\dd r
   \bigg\|_V
\\
 &
 \quad +
 \1_{\Omega_{R,t_{i-1} }}
   \bigg\|
   \int_0^s E(s-r) P_N  F(Y^{M,N}_{\lfloor r \rfloor}) 
   \tfrac{\tau \|P_N F(Y^{M,N}_{\lfloor r \rfloor})\|}{1 + \tau \| P_N F(Y^{M,N}_{\lfloor r \rfloor}) \|}
   \,\dd r
   \bigg\|_V
 +
   \| P_N \mathcal{O}_s  \|_V
\\
 & =:
   \| E_N( s ) Y_0^{M,N} \|_V
   +
   I_1 + I_2 
   + 
   \| P_N \mathcal{O}_s  \|_V.
\end{split}
\end{equation}
Before proceeding further, we claim
{\color{black}{
\begin{equation} \label{eq:Omega-Y-estimate}
\1_{\Omega_{R,t_{i-1} }} \| Y_r^{M,N} \|_V 
\leq 
C ( 
1
+
R 
+ 
\tau^{\frac34} R^{3} 
+ 
\| 
\mathcal{O}_r^N 
\|_V
+
\| 
\mathcal{O}_{ \lfloor r \rfloor }^N
\|_V 
),
\quad
 \forall \, r \in [0, t_i).
\end{equation}
By the definition of $Y_r^{M,N} $ and noting $\int_{ \lfloor r \rfloor }^r E ( r - u ) P_N \, \dd W_u 
=  \mathcal{O}_r^N - E ( r - \lfloor r \rfloor )  \mathcal{O}_{ \lfloor r \rfloor }^N$,  $r \in [0, t_i)$, one can write
\begin{equation}
Y_r^{M,N} 
=
E ( r - \lfloor r \rfloor ) Y_{\lfloor r \rfloor}^{M,N} 
+
\int_{ \lfloor r \rfloor }^r  \tfrac{ E ( r - u ) F_N ( Y_{ \lfloor u \rfloor }^{M,N} ) } 
{ 1 + \tau \| F_N ( Y_{ \lfloor u \rfloor }^{M,N} ) \| }  \, \dd u
+
\mathcal{O}_r^N - E ( r - \lfloor r \rfloor )  \mathcal{O}_{ \lfloor r \rfloor }^N.
\end{equation}
This together with the boundedness of the semigroup $E (t)$ in $V$, i.e.,
$\| E (t) \phi \|_V \leq \| \phi \|_V$ and \eqref{eq:lem-V-H-regularity} for $\gamma = 0$ 
promises
\begin{equation} \label{eq:Omega-Y-estimate1}
\begin{split}
\1_{\Omega_{R,t_{i-1} }} \| Y_r^{M,N} \|_V
& \leq
\1_{\Omega_{R,t_{i-1} }} 
\Big(
\| E ( r - \lfloor r \rfloor ) Y_{ \lfloor r \rfloor }^{M,N} \|_V
+
\int_{ \lfloor r \rfloor }^r \| E ( r - u ) F_N ( Y_{ \lfloor u \rfloor }^{M,N} ) \|_V \, \dd u
\\
& 
\quad
+
\| 
\mathcal{O}_r^N - E ( r - \lfloor r \rfloor )  \mathcal{O}_{ \lfloor r \rfloor }^N
\|_V
\Big)
\\
& \leq
\1_{\Omega_{R,t_{i-1} }} 
\Big(
\| Y_{ \lfloor r \rfloor }^{M,N} \|_V
+
\int_{ \lfloor r \rfloor }^r ( r - u )^{-\frac14} \|  F ( Y_{ \lfloor u \rfloor }^{M,N} ) \| \, \dd u
+
\| 
\mathcal{O}_r^N 
\|_V
+
\| 
\mathcal{O}_{ \lfloor r \rfloor }^N
\|_V
\Big)
\\
& \leq
C 
\big( 
1
+
R 
+ 
\tfrac43 \tau^{3/4} R^{3} 
+ 
\| 
\mathcal{O}_r^N 
\|_V
+
\| 
\mathcal{O}_{ \lfloor r \rfloor }^N
\|_V 
\big)
.
\end{split}
\end{equation}
%
This shows the claim \eqref{eq:Omega-Y-estimate}.
}}
With the aid of \eqref{eq:F-one-sided-condition} and \eqref{eq:Omega-Y-estimate},  the first term $I_1$ can be treated as follows,
\begin{equation} \label{eq:I1-first}
\begin{split}
   I_1
 &
   \leq \1_{\Omega_{R,t_{i-1} }}
   \int_0^s (s-r)^{- \frac{1}{4}}
   \| F(Y^{M,N}_r)
   -
   F(Y^{M,N}_{\lfloor r \rfloor}) \|
   \,\dd r
\\
 &
   \leq \1_{\Omega_{R,t_{i-1} }} C 
   \int_0^s 
   (s-r)^{- \frac{1}{4}} 
   ( 1 + R^2 + \tau^\frac32 R^6 + {\color{black}{ \| 
\mathcal{O}_r^N 
\|_V^2
+
\| 
\mathcal{O}_{ \lfloor r \rfloor }^N
\|_V^2 }})
   \| Y^{M,N}_r - Y^{M,N}_{\lfloor r \rfloor} \|
   \,\dd r,
\end{split}
\end{equation}
where $r \in [0, s]$, $ s \in [0, t_i]$,
\begin{equation} 
\begin{split}
   Y^{M,N}_r - Y^{M,N}_{\lfloor r \rfloor}
 & =
   [ E ( r ) - E ( \lfloor r \rfloor ) ] Y^{ M, N}_0 
   +
   \int_0^r E(r-u) \tfrac{ P_N F (Y^{M,N}_{\lfloor u \rfloor}) } { 1 + \tau \| P_N F (Y^{M,N}_{\lfloor u \rfloor})  \|  }
   \,\dd u
\\
 &
 \quad 
   -
   \int_0^{\lfloor r \rfloor} E({\lfloor r \rfloor}-u) \tfrac{ P_N F (Y^{M,N}_{\lfloor u \rfloor}) } { 1 + \tau \| P_N F (Y^{M,N}_{\lfloor u \rfloor})  \|  }
   \,\dd u
   + 
   P_N \mathcal{O}_r
   -
   P_N \mathcal{O}_{\lfloor r \rfloor}.
\end{split}
\end{equation}
%
This ensures that
\begin{equation}
\begin{split}
 & \1_{\Omega_{R,t_{i-1} }} \| Y^{M,N}_r - Y^{M,N}_{\lfloor r \rfloor} \|
   \\
 & \quad
   \leq \tau^{\frac{ \beta }{2}} \| Y^{M,N}_0 \|_\beta 
   + 
   1_{\Omega_{R,t_{i-1}}} \Big\| \int_0^{\lfloor r \rfloor} E({\lfloor r \rfloor}-u)
   \big(
   E( r - {\lfloor r \rfloor})-I \big) \tfrac{ P_N F (Y^{M,N}_{\lfloor u \rfloor}) } { 1 + \tau \| P_N F (Y^{M,N}_{\lfloor u \rfloor})  \|  }
   \,\dd u
   \Big\|
\\
 &
 \qquad +
   \1_{\Omega_{R,t_{i-1} }}
   \Big\| \int_{\lfloor r \rfloor}^r E ( r - u ) \tfrac{ P_N F (Y^{M,N}_{\lfloor u \rfloor}) } { 1 + \tau \| P_N F (Y^{M,N}_{\lfloor u \rfloor})  \|  }
   \,\dd u
   \Big \|
   +
   \1_{\Omega_{R,t_{i-1} }} \| P_N ( \mathcal{O}_r - \mathcal{O}_{\lfloor r \rfloor} ) \|
\\
 &
    \leq \tau^{\frac{ \beta }{ 2 } } \| X_0 \|_\beta 
    + 
    C(1 + R^ { 3 } ) (\tau^{ \frac34 } + \tau)
 +
   \| \mathcal{O}_r - \mathcal{O}_{\lfloor r \rfloor} \| .
\end{split}
\end{equation}
Inserting this into \eqref{eq:I1-first} results in
\begin{equation}
\begin{split}
   I_1
 &
   \leq C \int_0^s 
   ( 1 + R^2 + \tau^\frac32 R^6 + {\color{black}{ \| 
\mathcal{O}_r^N 
\|_V^2
+
\| 
\mathcal{O}_{ \lfloor r \rfloor }^N
\|_V^2 }}
) 
   ( s - r )^{- \frac{1}{4} }
   \Big [
   \tau^{\frac{\beta}{2}} \| X_0 \|_\beta
   + C( 1+R^{ 3 } ) \tau^{ \frac34 }
   + \| \mathcal{O}_r - \mathcal{O}_{\lfloor r \rfloor} \|
   \Big]
   \, \dd r
\\
 &
   = C \tau^{\frac{\beta}{2}} \| X_0 \|_\beta 
   \int_0^s 
   ( 1 + R^2 + \tau^\frac32 R^6 + {\color{black}{ \| 
\mathcal{O}_r^N 
\|_V^2
+
\| 
\mathcal{O}_{ \lfloor r \rfloor }^N
\|_V^2 }}
)
( s - r )^ {-\frac14} \, \dd r 
 \\
 &
 \quad  
 +
   C (1 + R^{3 } ) \tau^{ \frac34 } 
   \int_0^s
   ( 1 + R^{ 2 } + \tau^{\frac32} R^6 + {\color{black}{ \| 
\mathcal{O}_r^N 
\|_V^2
+
\| 
\mathcal{O}_{ \lfloor r \rfloor }^N
\|_V^2 }} 
) ( s - r )^ {-\frac14} \, \dd r 
\\
 &
 \quad 
 +
   C \! \int_0^s ( s - r)^{- \frac14} ( 1 + R^2 + \tau^\frac32 R^6 
   +
   {\color{black}{ \| 
\mathcal{O}_r^N 
\|_V^2
+
\| 
\mathcal{O}_{ \lfloor r \rfloor }^N
\|_V^2 }}
) \| \mathcal{O}_r - \mathcal{O}_{\lfloor r \rfloor} \|
   \,\dd r.
\end{split}
\end{equation}
Therefore, letting $ R = R^{\tau} := \tau^{-\frac{\beta}{ 4 }} $ 
and considering \eqref{eq:lem-stoch-conv-regularity}-\eqref{eq:lem-stoch-conv-spatial-regularity} 
one can further infer that
\begin{equation}\label{eq:I1-moment-bound}
\| I_1 \|_{L^{9p} (\Omega, \R) } 
\leq 
C 
( 
1 
+ 
\| X_0 \|_{L^{9p} (\Omega, \dot{H}^{\beta} ) }
+
\|A^{\frac{\beta-1}{2}}  \|_{\mathcal{L}_2 (H) } ).
\end{equation}
%
%
%
In a similar manner, choosing $ R = R^{\tau} := \tau^{-\frac{\beta}{ 4 }} $ enables us to treat $I_2$ as follows:
\begin{equation} \label{eq:I2-estimate}
\begin{split}
   I_2
 &
   \leq \1_{\Omega_{R,t_{i-1}}} \int_0^s \| E( s -r ) P_N F(Y^{M,N}_{\lfloor r \rfloor}) \|_V \cdot \tau \| P_N F(Y^{M,N}_{\lfloor r \rfloor})\|
   \,\dd r
\\
 &
   \leq \1_{\Omega_{R,t_{i-1}}} \tau \int_0^s ( s - r )^{ -\frac{ 1 }{ 4 } } \| F(Y^{M,N}_{\lfloor r \rfloor}) \|^2
   \,\dd r
\\
 &   
\leq C (R^{6} + 1) \tau
\\
 &
 \leq C ( \tau^{ \frac{ 2 - 3\beta}{2}  } + \tau ),
\end{split}
\end{equation}
where we also used \eqref{eq:lem-V-H-regularity} with $\gamma = 0$ 
{\color{black}{
and the structure of the taming scheme, namely,
$\tfrac{ z }{ 1 + \tau z } \leq z,  z \geq 0$.}}
%
Bearing \eqref{eq:I1-moment-bound}, \eqref{eq:I2-estimate} and \eqref{eq:lemma-PN-Stoch-Conv-V-bound} in mind, 
one can deduce from \eqref{eq:Z-V-subevents} that, for any $ s \in [0, t_i]$,
\begin{equation} \label{eq:Z-MB-subevent}
   \E [ \1_{\Omega_{R,t_{i-1} }}
   \| Z^{M,N}_s \|_V ^{ 9p} ] \leq C < \infty.
\end{equation}
This together with \eqref{eq:Zt-Yt} and \eqref{eq:Z-MB-subevent} immediately implies
\begin{equation}
   \E [ \1_{\Omega_{R^{\tau},t_{i-1}}} \| \bar{Y}^{M,N}_{t_i} \|_V^p] \leq C < \infty .
\end{equation}
Combining this with \eqref{eq:defn-Ybar} verifies the desired assertion \eqref{eq:lem-moment-bound0}.
$\square$

{\color{black}{As a direct consequence of Lemma \ref{Lem:Full-discrete-MB}, the following corollary
offers a priori moment bounds of the numerical approximation on subevents $ \Omega_{ R^{ \tau },t_i}, i \in \{ 0, 1,..., M \} $.
\begin{cor}
\label{cor:moment-bounds}
Let $ p \in [2, \infty)$ and $ R^{\tau} := \tau^{-\frac{\beta}{ 4 }} $ for any $ \beta \in (0, \tfrac12) $. 
Let Assumptions \ref{ass:A}-\ref{ass:X0}, \ref{ass:X0-additional} hiold
and let the approximation process $Y^{M,N}_{t_i}, i \in \{ 0, 1,..., M \}$  be produced by \eqref{eq:full.Tamed-AEE}.
Then it holds
\begin{equation}
\label{eq:cor-moment-bound}
\sup_{ M, N \in \N} \sup_{ i \in \{ 0, 1,..., M \}} 
\E\big[ \1_{\Omega_{ R^{ \tau }, t_{i} } } \| Y^{M,N}_{t_i} \|_V^{p} \big]
<
\infty.
\end{equation}
\end{cor}
{\it Proof of Corollary \ref{cor:moment-bounds}. }
Thanks to the assertion \eqref{eq:lem-moment-bound0} and the fact that $ \Omega_{R,t_i} \subset  \Omega_{R,t_{i-1}} $, one can easily deduce
\begin{equation} \label{eq:estimate-large-events}
\sup_{ M, N \in \N} \sup_{ i \in \{ 0, 1,..., M \}} 
\E\big[ \1_{\Omega_{ R^{ \tau }, t_{i} } } \| Y^{M,N}_{t_i} \|_V^{p} \big]
\leq
\sup_{ M, N \in \N} \sup_{ i \in \{ 0, 1,..., M \}} 
\E\big[ \1_{\Omega_{ R^{ \tau }, t_{i-1} } } \| Y^{M,N}_{t_i} \|_V^{p} \big]
<
\infty,
\end{equation}
as required. $\square$

Equipped with the bounded moments \eqref{eq:cor-moment-bound} on subevents $ \Omega_{ R^{ \tau },t_i}, i \in \{ 0, 1,..., M \} $, 
one only needs to deduce bounded moments on its complement $\Omega_{R^{\tau},t_m}^c$
before arriving at the expected bounded moments on the whole probability space stated as follows.
}}
\begin{thm}\label{thm:Numer-Moment-Bound}
Let Assumptions \ref{ass:A}-\ref{ass:X0}, \ref{ass:X0-additional} be fulfilled
and let the approximation process $Y^{M,N}_{t_i}, i \in \{ 0, 1,..., M \}$  be produced by \eqref{eq:full.Tamed-AEE}. 
Then for any $p \in [2, \infty)$ it holds that
\begin{equation} \label{eq:thm-Numer-Moment-Bound}
\sup_{ M, N \in \N} \sup_{ m \in \{ 0, 1,..., M \}} 
\E\big[ \| Y^{M,N}_{t_m} \|_V^{p} \big]
<
\infty.
\end{equation}
\end{thm}
{\it Proof of Theorem \ref{thm:Numer-Moment-Bound}.}
As discussed above, it remains to bound
$
\sup_{ M, N \in \N} \sup_{ m \in \{ 0, 1,..., M \}} 
\E [ \1_{\Omega_{R^{\tau},t_m}^c} \| Y^{M,N}_{t_m} \|_V^p]. 
$
{\color{black}{Owing to the boundedness of the semigroup $E (t)$ in $V$, i.e.,
$\| E (t) \phi \|_V \leq \| \phi \|_V$, \eqref{eq:lem-V-H-regularity} with $\gamma = 0$
and the structure of the taming scheme, i.e. $\tfrac{ z }{ 1 + \tau z } \leq \tau^{-1}$ for $z \geq 0$,}} 
one can infer 
\begin{equation} \label{eq:numer-bound-1}
\begin{split}
   \|Y_{t_m}^{M,N}\|_V
 &
   \leq \|E( t_m ) P_N X_0 \|_V
   +
   \| P_N \mathcal{O}_{t_m} \|_V
   +
   \int_0^{t_m}
  \Big \|
   E( t_m - s ) \tfrac{P_N \, F(Y^{M,N}_{\lfloor s \rfloor})}{1+\tau\| P_N \, F(Y^{M,N}_{\lfloor s \rfloor}) \|}
  \Big \|_V
   \,\dd s
\\
 &
 {\color{black}{
   \leq  \| P_N X_0 \|_V
   + 
   \int_0^{t_m}
   (t_m - s)^{-\frac 1 4}
   \tfrac{  \| P_N
   \, F(Y^{M,N}_{\lfloor s \rfloor}) \| }
   {1 + \tau\| P_N F(Y^{M,N}_{\lfloor s \rfloor}) \|}
   \,\dd s
   + \| P_N \mathcal{O}_{t_m} \|_V
   }}
\\
 &
   \leq 
   \| P_N X_0 \|_V
   + \tfrac 43 t_m^{\frac34} \tau^{-1}
   + 
   {\color{black}{\| P_N \mathcal{O}_{t_m} \|_V}},
   \quad
   m \in \{ 0, 1, 2,..., M \}.
\end{split}
\end{equation}
{\color{black}{
Thanks to Lemma \ref{Lem:PN-Stoch-Conv-V-bound} and Assumption \ref{ass:X0-additional}, we can show that for any $p \geq 2$,
\begin{equation} \label{eq:Yt-crude-estimate}
\|Y_{t_m}^{M,N}\|_{L^p( \Omega, V)}
\leq
C ( 1 + \tau^{-1} ),
\quad
   m \in \{ 0, 1, 2,..., M \}.
\end{equation}
}}
%
Meanwhile, one can learn that
\begin{equation}
   \Omega^c_{R^{\tau},t_m}
   =
   \Omega^c_{R^{\tau},t_{m-1}}
   +
   \Omega_{R^{\tau},t_{m-1}}
   \cdot
   \{\omega\in\Omega : \| Y^{M,N}_{t_m}\|_V > R^{\tau } \},
\end{equation}
and as a result
\begin{equation}
\begin{split}
   \1_{\Omega^c_{R^{\tau},t_m}}
& =
   \1_{\Omega^c_{R^{\tau},t_{m-1}}}
   +
   \1_{\Omega_{R^{\tau},t_{m - 1} } }
   \cdot
   \1_{\{\| Y^{M,N}_{t_m}\|_V > R^{\tau}\}}
=
   \sum_{i=0}^m \1_{\Omega_{R^{\tau}, t_{i-1} } }
   \cdot
   \1_{\{\| Y^{M,N}_{t_i}\|_V > R^{\tau}\}} ,
\end{split}
\end{equation}
where we recall
$ 
\1_{\Omega^c_{R^{\tau},t_{-1} }} = 0.
$
By \eqref{eq:Yt-crude-estimate} and the Chebyshev inequality, one can show, 
for any $m \in \{ 0, 1,...,M \}$, $M \in \N$,
\begin{equation} \label{eq:moment-bound-part2}
\begin{split}
   \E
   \big [
   \1_{\Omega^c_{R^{\tau},t_m}} \|Y_{t_m}^{M,N} \|_V^p
   \big ]
 &
 =
   \sum_{i=0}^m \E
   \big [
   \|Y_{t_m}^{M,N} \|_V^p \cdot
   \1_{\Omega_{ R^{\tau},t_{i-1} } }
   \1_{\{\| Y^{M,N}_{t_i}\|_V > R^{\tau}\}}
   \big ]
 \\
 &
   \leq \sum_{i=0}^m
   \Big ( \E
   \big [
   \|Y_{t_m}^{M,N} \|_V^{2p}
   \big ]
   \Big )^{\frac{1}{2}}
   \Big(
   \E
   \big [
   \1_{\Omega_{R^{\tau}, t_{i-1} }}
   \1_{\{\| Y^{M,N}_{t_i}\|_V > R^{\tau}\}}
   \big ]
   \Big )^{\frac{1}{2}}
 \\
 &
 {\color{black}{
   \leq \sum_{i=0}^m C ( 1+\tau^{-p})
   \Big(
   \mathbb{P} \big(
   \omega \in \Omega_{R^{\tau},t_{i-1} } : \| Y^{M,N}_{t_i}\|_V
   > R^{\tau}
   \big)
   \Big)
   ^{\frac{1}{2}}
   }}
 \\
 &
 {\color{black}{
   = 
   C ( 1+\tau^{-p})
   \sum_{i=0}^m
   \Big(
   \mathbb{P} \big(
   \omega \in \Omega :
   \1_{\Omega_{R^{\tau},t_{i-1} }} \| Y^{M,N}_{t_i}\|_V
   > R^{\tau}
   \big)
   \Big)
   ^{\frac{1}{2}}
   }}
\\
 &
   \leq C ( 1+\tau^{-p})
   \sum_{i=0}^m
   \Big (
   \E
   \Big [
   \1_{\Omega_{R^{\tau}, t_{i-1} }} \| Y^{M,N}_{t_i}\|_V^{ \frac{8(p+1)}{\beta}} / (R^{\tau})^{ \frac{8(p+1)}{\beta} }
   \Big ]
   \Big )
   ^{\frac{1}{2}}
\\
 &
   \leq
    C ( 1+\tau^{-p})
   \sum_{i=0}^m
   \tau^{ p+1 }
   \Big (
   \E
   \Big [
   \1_{\Omega_{R^{\tau}, t_{i-1} }} \| Y^{M,N}_{t_i}\|_V^{ \frac{8(p+1)}{\beta}}
   \Big ]
   \Big )
   ^{\frac{1}{2}}
   < \infty.
\end{split}
\end{equation}
This estimate together with \eqref{eq:cor-moment-bound} yields
the required estimate \eqref{eq:thm-Numer-Moment-Bound}.
$\square$

With Theorem \ref{thm:Numer-Moment-Bound} at hand, one can use standard arguments to obtain the coming corollaries.
%
%
\begin{cor} \label{cor:numer-spatial-regularity}
Let Assumptions \ref{ass:A}-\ref{ass:X0}, \ref{ass:X0-additional} be fulfilled
and let the approximation process $Y^{M,N}_{t}, t \in [0, T]$ be given by \eqref{eq:full-discrete-continous}.
Then for any $p \in [2, \infty)$ and $\beta < \tfrac12$ we have
\begin{equation}
\sup_{ M, N \in \N, \, t \in [0, T]}
\| 
Y^{M,N}_{t}
\|
_{ L^p ( \Omega, \dot{H}^{\beta} ) }
+
\sup_{ M, N \in \N, \, t \in [0, T]}
\| 
Y^{M,N}_{t}
\|
_{ L^p ( \Omega, V ) }
<
\infty.
\end{equation}
\end{cor}
\begin{cor}\label{cor:numerical-holder-regularity}
Let Assumptions \ref{ass:A}-\ref{ass:X0}, \ref{ass:X0-additional} be fulfilled
and let the approximation process $Y^{M,N}_{t}, t \in [0, T]$ be given by \eqref{eq:full-discrete-continous}. 
Then for any $p \in [2, \infty)$ and $\beta < \tfrac12$ {\color{black}{there exists a constant $C$, 
depending on $p, \beta, T$ and the initial data but not depending on $M,N$, such that}}
\begin{equation}\label{eq:numer-regularity-normal}
      \| Y_t^{M,N} - Y^{M,N}_s \|_{ L^p ( \Omega, H ) }
     \leq C (t-s)^{ \frac{\beta}{2} } , 
     \quad
     0 \leq s<t \leq T.
\end{equation}
\end{cor}
%
%
{\color{black}{\subsection{Refined temporal H\"{o}lder regularity results in negative Sobolev spaces}
\label{subset:negative-sobleve}
In addition to the a priori moment bounds for the approximations, we further rely on 
refined temporal H\"{o}lder regularity results in negative Sobolev spaces, as stated in Corollary  \ref{cor:F-diff-regularity},
which are essentially used in proving the temporal convergence rates of order almost $\tfrac12$.
To arrive at Corollary \ref{cor:F-diff-regularity}, we also need commutativity properties of the nonlinearity 
 and the improved temporal H\"{o}lder regularity of the solution in negative Sobolev spaces, 
as presented in Lemma \ref{Lem:F-negative-soblev} and Lemma \ref{Lem:regularity-negative-soblev}, respectively.

For the first step, we reveal commutativity properties of the nonlinearity in the forthcoming lemma.
}}
\begin{lem}\label{Lem:F-negative-soblev}
Let 
$F \colon L^{6} ( D; \R) \rightarrow H$ 
be a mapping determined by Assumption \ref{ass:F}.
Then for any $\beta \in (0, \tfrac12)$ and $ \eta > \tfrac12$ {\color{black}{there exists a constant $C$
depending on $\beta, \eta, \{a_i\}_{i=0}^{3}$ such that}}
\begin{equation} \label{eq:lem-derivative-F-regularity}
\| F'( \chi ) \nu \|_{ - \eta }
\leq
C \big(1+\max{ \{ \| \chi \|_V ,\| \chi \|_{\beta} \} }
   ^{2} \big)  \| \nu \|_{-\beta},
\quad
\forall
\chi \in V \cap \dot{H}^{\beta}, \nu \in V.
\end{equation}
\end{lem}
{\it Proof of Lemma \ref{Lem:F-negative-soblev}. }
As $ \beta \in (0, \frac{1}{2} ) $, standard arguments with the Sobolev-Slobodeckij norm (see, e.g., \cite[(19.14)]{thomee2006galerkin}) 
and properties of the nonlinear mapping guarantee
\begin{equation}
\begin{split}
   \| F' ( \chi ) \psi \|_{\beta}^2
 & \leq
   C  \|F' ( \chi ) \psi \|^2
   +
   C \int_0^1 \int_0^1 \frac{\big|f'(\chi(x)) \psi(x)
   -
   f'(\chi(y)) \psi(y)\big|^2} {|x-y|^{2{\beta}+1}} \, \dd y \dd x
\\
 &
   \leq
   C  \|F ' ( \chi ) \psi \|^2
   +
   C \int_0^1 \int_0^1 \frac{\big|f'(\chi(x)) (\psi(x)-\psi(y))\big|^2}
   {|x-y|^{2{\beta}+1}}
   \, \dd y \dd x
\\
 &
   \quad +
   C \int_0^1 \int_0^1
   \frac{ \big| [ f'(\chi(x))- f'(\chi(y)) ] \psi(y)\big|^2}
   {|x-y|^{2{\beta}+1}}
   \, \dd y \dd x
\\
 &
  \leq
   C  \big\|F ' ( \chi ) \psi \big\|^2
   +
   C  \big\| f'(\chi(\cdot)) \big\|_V^2 \cdot \| \psi \|^2_{W^{{\beta},2}}
 +
   C
   \big\|
   f''( \chi(\cdot) )
   \big\|_V^2
   \cdot
   \|\psi\|_V^2
   \cdot
   \| \chi \|^2_{W^{{\beta},2}}
\\
 &
   \leq
   C \big( 1+ \|\chi\|_V^{ 4 } \big)  \|\psi\|^2
   +
   C  \big( 1+ \|\chi\|_V^{ 4 } \big) \|\psi\|^2_{\beta}
 +
   C \big( 1 +  \|\chi\|_V^{ 2 } \big)
   \|\psi\|^2_V \cdot \|\chi\|^2_{\beta}
\\
 &
   \leq
   C  \big(1+\max{ \{ \| \chi \|_V ,\| \chi \|_{\beta} \} }
   ^{4} \big)
   (\|\psi\|^2_{\beta} + \|\psi\|^2_V ).
\end{split}
\end{equation}
Accordingly, for any $\beta \in (0, \tfrac12)$ and $ \eta > \tfrac12$, one uses the Sobolev embedding inequality to derive
\begin{equation}
\begin{split}
   \| F' ( \chi ) \nu \|_{ - \eta}
 &
   =
   \sup_{\| \varphi  \| \leq 1 }  \Big|\big\langle A^{-\frac{\eta}{2}}
   F' ( \chi ) \nu ,\varphi \big\rangle  \Big|
   =
   \sup_{\| \varphi  \| \leq 1 }  \Big|\big\langle \nu,
   ( F' ( \chi ) )^\ast \cdot A^{-\frac{\eta}{2}} \varphi \big\rangle \Big|
\\
 &
   =\sup_{\| \varphi  \| \leq 1 }  \Big|\big\langle A^{-\frac{\beta}{2}}
    \nu , A^{\frac{\beta}{2}} F' ( \chi )
    A^{-\frac{\eta}{2}}\varphi \big\rangle  \Big|
\\
 &
   \leq
   \sup_{\| \varphi  \| \leq 1 } \| \nu\|_{-\beta} \cdot
   \big\| F' ( \chi ) A^{-\frac{\eta}{2}}\varphi \big\|_{\beta}
\\
 &
   \leq
   \sup_{\| \varphi  \| \leq 1 } \| \nu\|_{-\beta} \cdot
   C \big(1+\max{ \{ \| \chi \|_V ,\| \chi \|_{\beta} \} }
   ^{ 2 } \big)
   (\| \varphi \|_{\beta-\eta}+ \| A^{-\frac{\eta}{2}}\varphi \|_V )
\\
 &
   \leq
   C_{\beta} \, \big(1+\max{ \{ \| \chi \|_V ,\| \chi \|_{\beta} \} }
   ^{ 2 } \big) \| \nu\|_{-\beta}.
\end{split}
\end{equation}
This completes the proof.
$\square$

{\color{black}{For the second step, we establish the improved temporal H\"{o}lder regularity of the solution in negative Sobolev spaces
as follows.}}
\begin{lem} \label{Lem:regularity-negative-soblev}
Let Assumptions \ref{ass:A}-\ref{ass:X0}, \ref{ass:X0-additional} be fulfilled
and let the approximation process $Y^{M,N}_{t}, t \in [0, T]$  be produced by \eqref{eq:full-discrete-continous}.
Then for any $p \in [2, \infty)$, {\color{black}{$\delta \in [0, \tfrac12]$}} and $\beta < \tfrac12$ 
{\color{black}{there exists a constant $C$ depending on $T, p, \delta, \beta, \{a_i\}_{i=0}^{3} $ 
and the initial data but independent of $M,N$ such that}}
\begin{equation}\label{eq:numer-regularity-negative}
      \| Y_t^{M,N} - Y^{M,N}_s \|_{ L^p ( \Omega, \dot{H}^{{\color{black}{ -\delta}} } ) }
     \leq C (t-s)^{ {\color{black}{ \frac{\beta + \delta }{2} }} } , 
     \quad
     0 \leq s < t \leq T.
\end{equation}
\end{lem}
{\it Proof of Lemma \ref{Lem:regularity-negative-soblev}}. 
The definition \eqref{eq:full-discrete-continous} implies, for $ 0 \leq s < t \leq T $,
\begin{equation} \label{eq:Y-diff-negative-sobleve0}
\begin{split}
   Y_t^{M,N}-Y^{M,N}_s
 &
   =
   \big(
   E_N ( t - s ) - I
   \big)
   Y^{M,N}_s
    \\ & \quad 
      +
   \int_s^t E_N ( t - r ) 
   \tfrac{ F_N(Y^{M,N}_{\lfloor r \rfloor})}{1+\tau
   \|
   F_N(Y^{M,N}_{\lfloor r \rfloor})
   \|}
   \, \dd r
   +
   \int_s^t E_N ( t - r ) P_N \,
   \dd W ( r ).
\end{split}
\end{equation}
Making use of \eqref{eq:E.Inequality}, \eqref{eq:AQ_condition} and the inequality 
$\| \Gamma \Gamma_1\|_{\mathcal{L}_2 } \leq \|\Gamma\|_{\mathcal{L} } \|\Gamma_1\|_{\mathcal{L}_2 }$, 
$  \Gamma \in \mathcal{L}(H), \Gamma_1 \in  \mathcal{L}_2(H) $ gives
\begin{equation} \label{eq:Y-diff-negative-sobleve1}
\begin{split}
\Big \|  \int_s^t E_N ( t - r ) P_N \,
   \dd W ( r ) 
 \Big \|_{ L^p ( \Omega, \dot{H}^{- \delta} ) }
 &
   \leq
   C
    \Big(
   \int_s^t \big \| A^{-\frac{\delta}{2}}
   E_N ( t - r ) P_N
   \big \|^2_{ \mathcal{L}_2(H) }
   \, \dd r
   \Big)^{\frac{ 1 }{2}}
   \leq
   C
   (t-s)^{ \frac{\beta + \delta }{2} }.
\end{split}
\end{equation}
A combination of \eqref{eq:E.Inequality} and Corollary \ref{cor:numer-spatial-regularity} shows
\begin{equation} \label{eq:Y-diff-negative-sobleve2}
\begin{split}
\| 
\big(
   E_N ( t - s ) - I
   \big)
   Y^{M,N}_s
\|_{ L^p ( \Omega, \dot{H}^{-\delta } ) }
 &
   \leq
   \big \|
   A^{- \frac{\beta + \delta }{2} }
   \big(
   E ( t - s ) - I
   \big) \big\|_{\mathcal{L} (H) }
   \cdot
   \|
   Y^{M,N}_s
   \|_{ L^p ( \Omega, \dot{H}^{\beta} ) }
\\
 &
   \leq
   C (t-s)^{ \frac{\beta + \delta }{2} }.
\end{split}
\end{equation}
Now we  proceed to estimate the remaining term in \eqref{eq:Y-diff-negative-sobleve0}, 
with the help of \eqref{eq:thm-Numer-Moment-Bound} and \eqref{eq:F-one-sided-condition},
\begin{equation} \label{eq:Y-diff-negative-sobleve3}
  \Big \| 
      \int_s^t E(t - r)  
   \tfrac{ F_N(Y^{M,N}_{\lfloor r \rfloor})}{1+\tau
   \|
   F_N(Y^{M,N}_{\lfloor r \rfloor})
   \|}
   \, \dd r
 \Big \|_{ L^p ( \Omega, \dot{H}^{-\delta} ) }
   \leq 
   C
   \int_s^t 
   \| F (Y^{M,N}_{\lfloor r \rfloor})
   \|_{ L^p ( \Omega, H ) }
   \, \dd r
   \leq
   C ( t - s ) .
\end{equation} 
Gathering \eqref{eq:Y-diff-negative-sobleve1}, \eqref{eq:Y-diff-negative-sobleve2} and \eqref{eq:Y-diff-negative-sobleve3} 
we deduce from \eqref{eq:Y-diff-negative-sobleve0} that \eqref{eq:numer-regularity-negative} is true.
$\square$

{\color{black}{
It is worthwhile to mention that  the approximate solution \eqref{eq:full-discrete-continous} enjoys higher order 
temporal H\"{o}lder regularity in negative Sobolev spaces $\dot{H}^\delta, \delta >0$,
as indicated by Lemma \ref{Lem:regularity-negative-soblev}. 
More precisely, the order of temporal H\"{o}lder regularity is increased to be $\tfrac{\beta + \delta}{2}$ when measured 
in the negative Sobolev space $\dot{H}^\delta, \delta \in (0, \tfrac12]$, in contrast to order $\tfrac{\beta}{2}$ in $H$ by 
Corollary \ref{cor:numerical-holder-regularity}. So taking $\delta = 0$ in \eqref{eq:numer-regularity-negative}
reduces into the usual temporal H\"{o}lder regularity results in $H$, i.e., \eqref{eq:numer-regularity-normal}
in Corollary \ref{cor:numerical-holder-regularity}.
%

Combining Lemmas \ref{Lem:F-negative-soblev}, \ref{Lem:regularity-negative-soblev} with Corollary \ref{cor:numer-spatial-regularity},
one can easily arrive at the following corollary, which is used several times in the error analysis,
to improve the strong convergence rate of the scheme.}}
\begin{cor} \label{cor:F-diff-regularity}
Let Assumptions \ref{ass:A}-\ref{ass:X0}, \ref{ass:X0-additional} be fulfilled
and let the approximation process $Y^{M,N}_{t}, t \in [0, T]$  be produced by \eqref{eq:full-discrete-continous}.
Then for any $p \in [2, \infty)$, $\beta < \tfrac12$ and $ \eta > \tfrac12 $,
{\color{black}{there exists a constant $C$ depending on $T, p, \eta, \beta, \{a_i\}_{i=0}^{3} $ 
and the initial data but independent of $M,N$ such that}}
\begin{equation}
\| F ( Y_t^{M,N} ) - F ( Y_s^{M,N} ) \|_{ L^p ( \Omega, \dot{H}^{ - \eta } ) }
\leq
C ( t - s )^{ \beta }, 
     \quad
     0<s<t<T.
\end{equation}
\end{cor}
\subsection{Main results: error bounds for the full discretization}
\label{subset:full-error-analysis}
%
{\color{black}{At the moment, we are well prepared to prove the expected strong convergence rate of the proposed scheme, 
as stated in the following theorem.
}}
\begin{thm}
[The space-time full error bounds]
\label{thm:full-discrete-scheme-error-bound}
Let Assumptions \ref{ass:A}-\ref{ass:X0}, \ref{ass:X0-additional} hold. 
Let $X (t) $ and $Y^{M,N}_{t}, t \in [0, T]$ be defined by 
\eqref{eq:mild-solution} and \eqref{eq:full-discrete-continous}, respectively.
For any $p \in [2, \infty)$ and $\beta \in (0, \tfrac12)$, 
{\color{black}{there exists a constant $C$ depending on $T, p, \beta, \{a_i\}_{i=0}^{3} $ 
and the initial data but independent of $M,N$ such that}}
\begin{equation} \label{eq:thm-main-result}
\sup_{ t \in [0, T]} 
\| X( t ) - Y^{M,N}_{t} \|_{L^p ( \Omega; H) } 
\leq
C \big(
     N ^{ - \beta } + \tau ^{  \beta }  \big).
\end{equation}
\end{thm}
%
%
{\it Proof of Theorem \ref{thm:full-discrete-scheme-error-bound}. }
Denoting $ e^{ M, N}_t : = P_N X(t) - Y_t^{ M, N } $, we note that
\begin{equation}\label{eq:main-thm-error-decom}
\| X( t ) - Y^{M,N}_{ t } \|_{L^p ( \Omega; H) } 
\leq
\| ( I - P_N ) X( t ) \|_{L^p ( \Omega; H) }
+
\| e^{M,N}_{t} \|_{L^p ( \Omega; H) },
\end{equation}
%
where 
\begin{equation}
   \frac{\dd }{\dd t} e^{ M, N}_t
   =
   - A_N
   e^{ M, N}_t
   +
   F_N
   \big (
   X ( t )
   \big )
   -
   \tfrac{ F_N(Y^{M,N}_{\lfloor t \rfloor})}{1+\tau
   \|
   F_N(Y^{M,N}_{\lfloor t \rfloor})
   \|} .
\end{equation}
This in conjunction with \eqref{eq:F-one-sided-condition} tells us that
\begin{equation} \label{eq:full-error-J0J1J2}
\begin{split}
   \big \|
   e^{ M, N}_t
   \big \|^p
   & =
   p \int_0^t
   \big \|
   e^{ M, N}_s
   \big \|^{p-2}
   \Big \langle
   e^{ M, N}_s, 
   - A_N
   e^{ M, N}_s
   +
   F_N
   \big (
   X (s) )
   -
   \tfrac{ F_N(Y^{M,N}_{\lfloor s \rfloor})}
   { 1 + \tau
   \|
    F_N ( Y^{M,N}_{\lfloor s \rfloor})
    \|}
    \Big \rangle
    \,\dd s
\\
 &
   \leq  
p \int_0^t
   \big \|
   e^{ M, N}_s
   \big \|^{p-2}
   \Big \langle
   e^{ M, N}_s, 
   F_N ( X (s) ) - F_N ( P_N X (s) ) 
   \\
 &
   \quad \quad \quad
   +
   F_N
   \big (
   Y^{M,N}_s
   \big )
    -
   \tfrac{ 
   F_N
   (
   Y^{M,N}_{\lfloor s \rfloor}
   )
   }
   { 1 + \tau
   \|
   F_N(Y^{M,N}_{\lfloor s \rfloor})
   \|}
   \Big \rangle
   \,\dd s
    +
    pL_0 \int_0^t
    \big \|
    e^{ M, N}_s
    \big \|^p
    \,\dd s
\\
 &  =
   p L_0 \int_0^t
   \big \|
    e^{ M, N}_s
    \big \|^p
    \,\dd s
   +
   p \int_0^t
   \big \|
   e^{ M, N}_s
   \big \|^{p-2}
   \big \langle
   e^{ M, N}_s, F( X (s) ) - F ( P_N X (s)  )
   \big \rangle
   \,\dd s
   \\
 &
   \quad +
   p \int_0^t
   \big \|
   e^{ M, N}_s
   \big \|^{p-2}
   \big \langle
   e^{ M, N}_s, F(Y_s^{M,N}) - F (Y^{M,N}_{\lfloor s \rfloor})
   \big \rangle
   \,\dd s
\\
 &
   \quad +
   p \int_0^t
   \big \|
   e^{ M, N}_s
   \big \|^{p-2}
   \Big \langle
   e^{ M, N}_s,
   \tfrac{ \tau
   \|
   F_N(Y^{M,N}_{\lfloor s \rfloor})
   \|
   \cdot
   F(Y^{M,N}_{\lfloor s \rfloor}) }{1+\tau
   \|
   F_N(Y^{M,N}_{\lfloor s \rfloor})
   \|}
   \Big \rangle
   \,\dd s
\\
 &  =:
   p L_0 \int_0^t
   \big \|
   e^{ M, N}_s
   \big \|^p
   \,\dd s
   +
   J_0
   + 
   J_1 
   + 
   J_2.
\end{split}
\end{equation}
Following the same lines as in estimates of \eqref{eq:spatial-error-before-final} and \eqref{eq:spatial-error-estimate-final}
one can bound the item $J_0$ as follows,
\begin{equation} \label{eq:estimateJ0}
\begin{split}
\E [ J_0 ] 
& \leq
p \, \E
\!
\int_0^t
   \|
   e^{ M, N}_s
   \|^{p-1}
   \|
   F( X (s) ) - F ( P_N X (s)  )
   \|
   \,\dd s
\\
& 
\leq
( p - 1 ) \int_0^t
   \E [
   \|
   e^{ M, N}_s
   \|^{p}
   ]
   \,
   \dd s
+
\int_0^t
 \E
 [
   \|
   F( X (s) ) - F ( P_N X (s)  )
   \|^p
 ]
 \, \dd s
\\
&
\leq
( p - 1 ) \int_0^t
   \E [
   \|
   e^{ M, N}_s
   \|^{p}
   ]
   \, \dd s
+
C
( \tfrac{1}{N} )^{ p \beta }.
\end{split}
\end{equation}
The term $J_2 $ is also easy to be treated, after taking 
the H\"{o}lder inequality and \eqref{eq:thm-Numer-Moment-Bound}
into account:
\begin{equation}\label{eq:estimateJ2}
\begin{split}
  \E [ J_2 ]
 &
   \leq p  \E
   \int_0^t
   \big \|
   e^{ M, N}_s
   \big \|^{p-1}
   \cdot
   \tau
   \big \|
   F(Y^{M,N}_{\lfloor s \rfloor})
   \big \|^2
   \,\dd s
\\
 &
   \leq
   (p-1)
   \int_0^t
   \E
   \big[
     \big \|
     e^{ M, N}_s
   \big \|^p
   \big]
   \,\dd s
   +
   \tau^p \int_0^t
   \E
   \big[
   \big \|
   F(Y^{M,N}_{\lfloor s \rfloor})
   \big \|^{2p}
   \big]
   \,\dd s
\\ &
   \leq
   (p-1) \int_0^t  \E
   \big [
   \big \|
   e^{ M, N}_s
   \big \|^p
   \big ]
   \,\dd s
   +
   C
   \tau^p .
\end{split}
\end{equation}
The remaining term $J_1$ must be handled more carefully.  
As usual, such a term is simply treated with the aid of temporal H\"{o}lder regularity of $ Y_s^{M,N} $
together with the Cauchy-Schwarz inequality and H\"{o}lder's inequality,
but to only attain order $\tau^{ \frac{\beta}{2} }$. 
In our analysis, we decompose $P_N X( s ) - Y_s^{ M, N }$ 
in the inner product into three parts.
To do so we recall that 
\begin{equation}
e^{ M, N}_s = P_N X ( s ) - Y_s^{ M, N } 
=
   \int_0^s E(s-r)
   \Big (
   F_N
   \big (
   X ( r )
   \big )
   -
   \tfrac{ F_N(Y^{M,N}_{\lfloor r \rfloor})}{1+\tau
   \|
   F_N(Y^{M,N}_{\lfloor r \rfloor})
   \|}
   \Big ) 
   \dd r,
\end{equation}
and split $J_1$ into three terms:
\begin{align}
\label{eq:J1-split}
   J_1
 & =
   p \int_0^t
   \big \|
   e^{ M, N}_s
   \big \|^{p-2}
   \Big \langle
   \int_0^s E(s-r)
   \Big (
   F_N
   \big (
   X ( r )
   \big )
   -
   \tfrac{ F_N(Y^{M,N}_{\lfloor r \rfloor})}{1+\tau
   \|
   F_N(Y^{M,N}_{\lfloor r \rfloor})
   \|}
   \Big )
   \,\dd r,
   F(Y_s^{M,N})-F(Y^{M,N}_{\lfloor s \rfloor})
   \Big \rangle
   \,\dd s
\nonumber \\
 & =
   p \int_0^t
   \big \|
   e^{ M, N}_s
   \big \|^{p-2}
   \Big \langle
   \int_0^s E(s-r)
   \big (
   F_N
   \big (
   X ( r )
   \big )
   -
   F_N(Y^{M,N}_{ r } )
   \big )
   \dd r ,
   F(Y_s^{M,N})-F(Y^{M,N}_{\lfloor s \rfloor})
   \Big \rangle
   \,\dd s
\nonumber \\
 &
   \quad +
   p \int_0^t
   \big \|
   e^{ M, N}_s
   \big \|^{p-2}
   \Big \langle
   \int_0^s E(s-r)
   \big (
   F_N(Y_r^{M,N})-F_N(Y^{M,N}_{\lfloor r \rfloor})
   \big)
   \,\dd r,
   F(Y_s^{M,N})-F(Y^{M,N}_{\lfloor s \rfloor})
   \Big \rangle
   \,\dd s
\nonumber \\
 &
   \quad +
   p \int_0^t
   \big \|
   e^{ M, N}_s \|^{p-2}
   \Big \langle
   \int_0^s E(s-r) \tfrac{ \tau
   \|
   F_N(Y^{M,N}_{\lfloor r \rfloor})
   \|
   \cdot
   F_N(Y^{M,N}_{\lfloor r \rfloor}) }{1+\tau
   \|
   F_N(Y^{M,N}_{\lfloor r \rfloor})
   \|}
   \,\dd r ,
   F(Y_s^{M,N})-F(Y^{M,N}_{\lfloor s \rfloor})
   \Big \rangle
   \,\dd s
\nonumber \\
 & =:
   J_{11} + J_{12} + J_{13}.
\end{align}
Since the estimates concerning $ J_{11} $ and $ J_{12} $ are demanding, we handle the item $ J_{13} $ first.
Utilizing \eqref{eq:F-one-sided-condition}, the H\"{o}lder inequality, \eqref{eq:thm-Numer-Moment-Bound} and  
\eqref{eq:numer-regularity-normal} results in
\begin{align}
{\color{black}{
   \E [ | J_{13} | ]
   }}
   &
   \leq p  \, \E \! \int_0^t \! \int_0^s
   \big \|
   e^{ M, N}_s
   \big \|^{p-2} 
   \cdot 
   \tau
   \big \|
   F(Y^{M,N}_{\lfloor r \rfloor})
   \big \|^2
   \cdot
   \big \| F(Y_s^{M,N})-F(Y^{M,N}_{\lfloor s \rfloor})
   \big \|
   \,\dd r \dd s
\nonumber \\
 & \leq
   C \tau \E \! \int_0^t \! \int_0^s
   \big \|
   e^{ M, N}_s
   \big \|^{p-2}
   \big \| F(Y^{M,N}_{\lfloor r \rfloor})
   \big \|^2 
   \big (
   1
   +
   \|
   Y_s^{M,N}
   \|_V^{ 2 }
   +
   \|
   Y^{M,N}_{\lfloor s \rfloor}
   \|_V^{ 2 }
   \big )
   \| Y_s^{M,N}-Y^{M,N}_{\lfloor s \rfloor} \|
   \,\dd r \dd s
\nonumber \\
 &
   \leq C \! \int_0^t \E 
   [
   \|
   e^{ M, N}_s
   \|^p
   ]
   \,\dd s
   +
   C \tau^{\frac{p}{2}} \int_0^t 
   \E
   \Big[
   \Big|
   \int_0^s
   \|
   F(Y^{M,N}_{\lfloor r \rfloor})
   \|^2
   \cdot
   \big (
   1+ \|Y_s^{M,N} \|_V^{2}
\nonumber \\
 &
   \quad +
   \| Y^{M,N}_{\lfloor s \rfloor}\|_V^{2}
   \big )
   \| Y_s^{M,N}-Y^{M,N}_{\lfloor s \rfloor} \|
   \,\dd r
   \Big |^{\frac{p}{2}}
   \Big]
   \,\dd s
\nonumber \\
 &
   \leq C \! \int_0^t \E
   [
   \|
   e^{ M, N}_s
   \|^p
   ]
   \,\dd s
   +
   C
   \tau^{\frac{p}{2}(1+  \frac{\beta}{2} )}.
\label{eq:EstimateJ13}
\end{align}
At the moment we come to the estimate of $ J_{11} $ and use the Taylor formula, 
the self-adjointness of operators $ F ' ( u )$ and $P_N$ to infer that
\begin{equation} \label{eq:J11-Estimate0}
\begin{split}
{\color{black}{
   \E [ | J_{11} | ]
   }}
 & =
   p \, \E
   \bigg [ \bigg|
   \int_0^t
   \big \|
   e^{ M, N}_s
   \big \|^{p-2}
   \Big \langle \! \int_0^s E(s - r)
   \big(
   F_N ( X ( r ) ) - F_N ( Y_r^{M,N} )
   \big )
   \,\dd r ,
   F(Y_s^{M,N})- F(Y^{M,N}_{\lfloor s\rfloor})
   \Big \rangle
   \,\dd s
   \bigg | \bigg]
\\
 & =
   p \, \E
   \bigg [ \bigg|
   \int_0^t
   \big \|
   e^{ M, N}_s
   \big \|^{p-2}
   \Big \langle \int_0^s E(s - r )
   P_N \! \int_0^1 F'
   \big( Y_r^{M,N} + \sigma ( X ( r ) - Y_r^{M,N} )
   \big)
   \dd \sigma
\\
 &
   \qquad   \qquad 
   \cdot
   \big(
   X ( r ) - Y_r^{M,N}
   \big)
   \,\dd r ,
   \
   F(Y_s^{M,N})- F(Y^{M,N}_{\lfloor s\rfloor})
   \Big \rangle
   \,\dd s
   \bigg | \bigg]
 \\
 & =
   p \, \E
   \bigg [ \bigg|
   \int_0^t \int_0^s \int_0^1
   \big \|
   e^{ M, N}_s
   \big \|^{p-2}
   \Big \langle
   X ( r ) - Y_r^{M,N},
   \big(
   F' \big( Y_r^{M,N} + \sigma ( X ( r ) - Y_r^{M,N} )
   \big)
   \big)^*
\\
 &
   \qquad   \qquad \qquad
    P_N E(s - r)
   \big[
   F(Y_s^{M,N})- F(Y^{M,N}_{\lfloor s\rfloor})
   \big]
   \Big \rangle
   \, \dd \sigma \,\dd r \, \dd s
   \bigg | \bigg]
\\
 &
   \leq
   p \, \E\int_0^t \int_0^s \int_0^1
   \big \|
   e^{ M, N}_s
   \big \|^{p-2}
   \cdot
   \big \|
   X ( r ) - Y_r^{M,N}
   \big \|
   \cdot \big \|
    F'
   \big( Y_r^{M,N} + \sigma ( X ( r ) - Y_r^{M,N} )
   \big)
\\
 &
   \qquad   \qquad \qquad   
    P_N A^{ \frac{\eta}{2} } E(s - r) A^{ - \frac{\eta}{2} }
   \big[
   F(Y_s^{M,N})- F(Y^{M,N}_{\lfloor s\rfloor})
   \big]
   \big \|
   \, \dd \sigma \,\dd r \, \dd s
\\
 &
   \leq
  C \, \E\int_0^t \int_0^s
   \big \|
   e^{ M, N}_s
   \big \|^{p-2}
   \big \|
   X ( r ) - Y_r^{M,N}
   \big \|
   \cdot 
   \big ( 1 + \| X ( r ) \|_V^2 + \| Y_r^{M,N} \|_V^2
    \big )
\\
 &
   \qquad   \qquad  \qquad   
\times ( s - r )^{- \frac{\eta}{2} }
   \big \|A^{-\frac{\eta}{2}}
   \big[
   F(Y_s^{M,N})- F(Y^{M,N}_{\lfloor s\rfloor})
   \big]
   \big \|
   \,\dd r \, \dd s.
\end{split}
\end{equation}
Further, employing {\color{black}{the Young inequalities $a b \leq \tfrac {a^\delta} {\delta} +  \tfrac {b^{\delta'} } {\delta'} $ for
$a , b \geq 0$, $\delta, \delta' > 1$, $\tfrac {1} {\delta} +  \tfrac {1 } {\delta'} = 1 $,}}, the H\"{o}lder inequality,
the transformation of integral domain and taking $\tfrac12 < \eta < 1$ give
\begin{equation}
\begin{split}
{\color{black}{
\E [ | J_{11} | ]
}}
 &
   \leq
   C \, \E \int_0^t \int_0^s
    ( s - r )^{ - \frac{ \eta }{ 2 } }
    \big \|
   e^{ M, N}_s
   \big \|^{p}
   \, \dd r \, \dd s
   +
   C \, \E 
   \int_0^t 
   \int_0^s 
   ( s - r )^{ - \frac{ \eta }{ 2 } }
   \big \|
   X ( r ) - Y_r^{M,N}
   \big \|^{\frac{p}{2}}
\\
 &
   \qquad   \qquad  \qquad 
   \times
\big ( 1 + \| X ( r ) \|_V^p + \| Y_r^{M,N} \|_V^p
    \big )
   \big \|A^{-\frac{\eta}{2}}
   \big[
   F(Y_s^{M,N})- F(Y^{M,N}_{\lfloor s\rfloor})
   \big]
   \big \|^{\frac{p}{2} }
   \,\dd r 
   \, \dd s
\\
 &
   \leq
   C \int_0^t \E
   \big [
     \big \|
   e^{ M, N}_s
     \big \|^{p}
   \big ]
   \,\dd s
   +
   C \, \E 
   \int_0^t 
   \int_0^s 
   ( s - r )^{ - \frac{ \eta }{ 2 } }
   \big \|
   X ( r ) - Y_r^{M,N}
   \big \|^{p}
   \, \dd r \, \dd s
   \\
   &
   \quad +
   C \, \E 
   \int_0^t 
   \int_0^s 
   ( s - r )^{ - \frac{ \eta }{ 2 } }
   \big ( 1 + \| X ( r ) \|_V^{2p} + \| Y_r^{M,N} \|_V^{2p}
    \big )
   \big \|A^{-\frac{\eta}{2}}
   \big[
   F(Y_s^{M,N})- F(Y^{M,N}_{\lfloor s\rfloor})
   \big]
   \big \|^{ p }
   \, \dd r \, \dd s
\\
 &
   \leq
   C \! \int_0^t \E
   \big [
     \big \|
   e^{ M, N}_s
     \big \|^{p}
   \big ]
   \,\dd s
   +
   C \! \int_0^t \E
   \big [
     \big \|
   X ( s )-Y_s^{M,N}
     \big \|^{p}
   \big ]
   \,\dd s
   \\ &
   \quad
   +
   C \!
   \int_0^t \int_0^s 
   ( s - r )^{ - \frac{ \eta }{ 2 } }
   \E
   \Big [ 
   \big ( 1 + \| X ( r ) \|_V^{2p} + \| Y_r^{M,N} \|_V^{2p}
    \big )
   \big \|A^{-\frac{\eta}{2}}
   \big[
   F(Y_s^{M,N})- F(Y^{M,N}_{\lfloor s\rfloor})
   \big]
   \big \|^p
   \Big ]
    \dd r \,\dd s,
\end{split}
\end{equation}
{\color{black}{where the first step is trivial for $p = 2$ due to \eqref{eq:J11-Estimate0} and holds true for $p > 2$ 
after using the aforementioned Young inequality with $\delta = \frac{p}{p-2} $, $\delta' = \frac{p}{2}$.
}}
To proceed further,
we resort to Corollary \ref{cor:F-diff-regularity} as well as \eqref{eq:optimal.regularity1},  
\eqref{eq:thm-Numer-Moment-Bound} and achieve
\begin{equation}\label{eq:EstimaeJ11}
\begin{split}
{\color{black}{
   \E[ | J_{11} | ]
   }}
 &
   \leq
   C  \! \int_0^t \E \big [
   \big\| e^{ M, N}_s \big \|^p \big] \, \dd s
   +
   C \lambda_{N + 1}^{- \frac{ p \beta }{2} }
   +
   C  \! \int_0^t \!
   \Big( 
   \E 
   \Big[
    \big \|A^{-\frac{\eta}{2}}
   \big[
   F(Y_s^{M,N})- F(Y^{M,N}_{\lfloor s\rfloor})
   \big]
   \big \|^{2p} 
   \Big]
   \Big)^{\frac{1}{2}}\, \dd s
\\
 &
   \leq
   C \! \int_0^t \E \, \big [ 
   \big\| e^{ M, N}_s \big \|^p \big ] \, \dd s
   +
   C \, ( \tfrac{1}{N} )^{ p \beta } 
   +
   C \, \tau^{p \beta} .
\end{split}
\end{equation}
%
Finally, it remains to deal with the estimate of  $ J_{12} $.
By {\color{black}{the Young inequalities $a b \leq \tfrac {a^\delta} {\delta} +  \tfrac {b^{\delta'} } {\delta'} $ for
$a , b \geq 0$, $\delta, \delta' > 1$, $\tfrac {1} {\delta} +  \tfrac {1 } {\delta'} = 1 $}}
and putting $\tfrac12 < \eta < 1$ one can derive that
\begin{align}
{\color{black}{
   \E[ | J_{12} | ]
   }}
 &
   =
   p\E  
   \bigg [ \bigg|
   \int_0^t   
   \|
   e^{ M, N}_s  
   \| ^ {p-2}
   \Big  \langle  \int_0^s  A^{\eta} E(s-r) A^{ -\frac{\eta}{2}}
   \Big ( F_N  \big ( Y^{M,N}_r  \big )
   -
   F_N  \big ( Y ^{M,N } _ { \lfloor r \rfloor }
    \big )  \Big )  \, \dd r,
\nonumber \\
 &
   \qquad  \qquad
   A ^ { -\frac{\eta}{2}}   \big (  F_N  \big ( Y^{M,N}_s  \big )
   -
    F_N \big ( Y^{M,N}_{ \lfloor s \rfloor }
    \big)  \big) \Big \rangle
    \, \dd s
    \bigg | \bigg]
\nonumber \\
 &
   \leq
    C \, \E  \int_0^t  \int_0^s
    (s-r)^{-\eta}
    \|
    e^{ M, N}_s  
    \|^{p-2}
    \cdot  \big \| A^{-\frac{\eta}{2}}
     \big(F  \big(  Y^{M,N}_r  \big)
     -
    F \big( Y^{M,N}_{ \lfloor r \rfloor }  \big) \big)  \big \|  
\nonumber \\
 &
  \qquad  \qquad
    \times
    \big \|  A^{-\frac{\eta}{2}} \big( F  \big(Y^{M,N}_s  \big)
    -
    F  \big( Y^{M,N}_{ \lfloor s \rfloor }  \big) \big) \big \|  
    \, \dd r \, \dd s
\nonumber \\
 &
   \leq
   C \int_0^t \int_0^s
   (s-r)^{-\eta}
   \E  [ \|
   e^{ M, N}_s
   \|^p  ]  \, \dd r \, \dd s
   +
   C \, \E  \int_0^t  \int_0^s (s-r)^{-\eta}  \big \| A^{-\frac{\eta}{2}}
   \big( F  \big( Y^{M,N}_r  \big)
   -
   F \big(Y^{M,N}_{ \lfloor r \rfloor } \big)\big)  \big \|^{\frac{p}{2}}
\nonumber \\
 &
   \qquad  \qquad \quad \, \,
   \times
   \big \| A^{-\frac{\eta}{2}} \big( F \big(Y^{M,N}_s  \big)
   -
   F \big( Y^{M,N}_{ \lfloor s \rfloor } \big) \big)  \big \|^{\frac{p}{2}}
    \, \dd r \, \dd s
\nonumber \\
  &
  \leq
  C \int_0^t 
   \E  [ \|
   e^{ M, N}_s
   \|^p  ]  \, \dd s
   +
   C  \int_0^t  \int_0^s (s-r)^{-\eta}  \E \big[  \big \| A^{-\frac{\eta}{2}}
   \big( F  \big( Y^{M,N}_r  \big)
   -
   F \big(Y^{M,N}_{ \lfloor r \rfloor } \big)\big)  \big \|^{p} \big ]
   \, \dd r \, \dd s
 \nonumber  \\ 
   & \quad
   +
   C  \int_0^t  \int_0^s (s-r)^{-\eta}  \E \big[  \big \| A^{-\frac{\eta}{2}}
   \big( F  \big( Y^{M,N}_s  \big)
   -
   F \big(Y^{M,N}_{ \lfloor s \rfloor } \big)\big)  \big \|^{p} \big ]
   \, \dd r \, \dd s
   ,
\end{align}
{\color{black}{where the last step but one is trivial for $p = 2$ and holds true for $p > 2$ after using the aforementioned 
Young inequality with $\delta = \frac{p}{p-2} $, $\delta' = \frac{p}{2}$. In the last step, we used the inequality
$a b \leq \tfrac{a^2}{2} + \tfrac{b^2}{2}, a, b \geq 0$. }}
Again, the use of Corollary \ref{cor:F-diff-regularity} leads us to
\begin{equation}\label{eq:EstimateJ12-final}
\begin{split}
{\color{black}{
\E[ | J_{12} | ] 
}}
    \leq 
C \int_0^t 
   \E  [ \|
   e^{ M, N}_s
   \|^p ]  \, \dd s
   +
   C \tau^{\beta p},
\end{split}
\end{equation}
which together with \eqref{eq:EstimateJ13}, \eqref{eq:EstimaeJ11}  forces us to recognize from \eqref{eq:J1-split} that
\begin{equation} 
{\color{black}{
\E[ | J_{1} | ] 
}}
\leq
C \int_0^t 
   \E  [ \|
   e^{ M, N}_s
   \|^p ]  \, \dd s
   +
   C \, ( \tfrac{1}{N} )^{ p \beta } 
   +
   C \tau^{\beta p}.
\end{equation}
Plugging this and \eqref{eq:estimateJ0}, \eqref{eq:estimateJ2}
into \eqref{eq:full-error-J0J1J2}, taking expectations and
applying the Gronwall inequality give
\begin{equation}
\| e^{M,N}_{t} \|_{L^p ( \Omega; H) }
\leq
C \big(
     N ^{ - \beta } + \tau ^{  \beta }  \big).
\end{equation}
Therefore, the desired error bound follows from \eqref{eq:main-thm-error-decom} by taking \eqref{eq:optimal.regularity2} 
and \eqref{eq:P-N-estimate} into account.
$\square$

%

%

%
%

%
%
%
 %
%
%
%
%
%
%
%
%
%
%
%
%
%
%
%
%
%
%
%
%
%
%
%

%

\section{Numerical experiments}
\label{sec:numerical-result}
Some numerical experiments are performed in this section to test the previous theoretical findings.  
Let us consider the 
stochastic Allen-Cahn equation with additive space-time white noise, described by
\begin{equation}\label{eq:num-result-AC-eqn}
\left\{
    \begin{array}{lll}
    \frac{\partial u}{\partial t} = \frac{\partial^2 u}{\partial x^2} + u - u^3 +  \dot{W} (t), \quad  t \in (0, 1], \:\: x \in (0,1),
    \\
     u(0, x) = \sin( \pi x), \quad\quad x \in (0,1),
     \\
     u(t, 0) = u(t,1) = 0,  \quad t \in (0, 1].
    \end{array}\right.
\end{equation}
Here $\{ W (t) \}_{t \in[0, T]}$ is a cylindrical $I$-Wiener process represented by 
\eqref{eq:Wiener-representation}.  In what follows, we will use the new fully discrete scheme 
\eqref{eq:full.Tamed-AEE} to approximate the continuous problem \eqref{eq:num-result-AC-eqn}.
Error bounds are always measured in terms of mean-square approximation errors at the endpoint $T = 1$, 
caused by spatial and temporal discretizations and the expectations are approximated by 
computing averages over 1000 samples.

Before proceeding further with simulations, it is helpful to mention that 
the stochastic convolution in the scheme \eqref{eq:num-result-AC-eqn} is easily implementable  
once one realizes that $\int_{t_m}^{t_{m+1}}\! E_N(t_{m+1}-s) P_N \mbox{d} W(s) = 
\sum_{i=1}^N \Lambda_i e_i$, where $  \Lambda_i = \int_{t_m}^{t_{m+1}} e^ {  - (t_{m+1}-s) \lambda_i } \dd \beta(s), 
1 \leq i \leq N $ are independent, zero-mean normally distributed random variables with explicit variances
$\E [ | \Lambda_i |^2 ] = \tfrac{ 1 - e^{ - 2 \lambda_i \tau } } { 2 \lambda_i }$. For more details on the implementation
of so-called AEE schemes, one can consult \cite[section 3]{jentzen2009overcoming} and \cite[section 4.1]{wang2014higher}.
%

\begin{figure}[htp]
\centering
      \includegraphics[width=4in,height=3in] {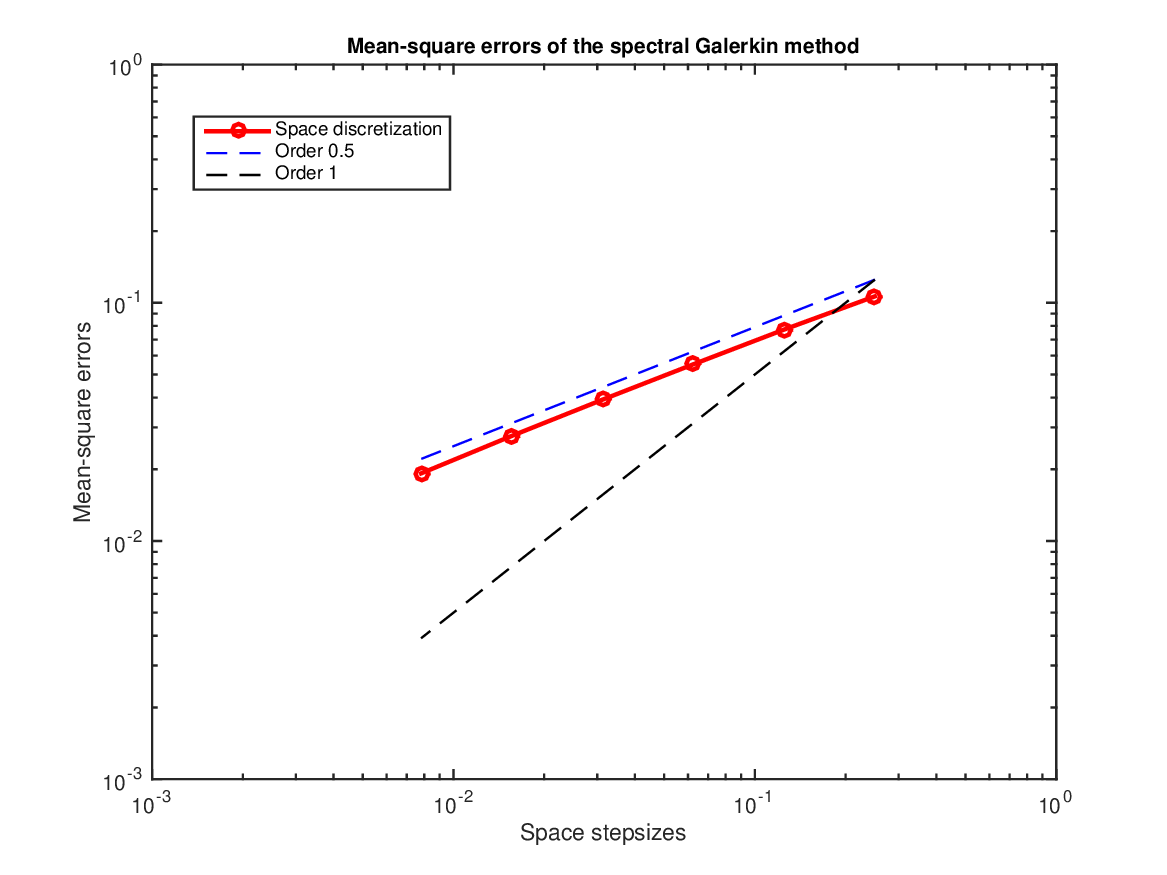}
 \caption{The convergence rate of the spectral Galerkin spatial discretization.}
\label{fig:spatial-error}
\end{figure}

To visually inspect the convergence rates in space, we identify the ``exact'' solution by using the full 
discretization with $M_{\text{exact}} = N_{\text{exact}} = 2^{11} = 2048$. 
The spatial approximation errors $\| X(1) - X^N(1) \|_{L^2 ( \Omega; H) }$ with
$ N= 2^i, i = 2, ..., 7 $ are depicted in Fig.\ref{fig:spatial-error}, against $\tfrac{1}{N}$ on a log-log scale, where one can observe
that the resulting spatial errors decrease at a slope close to $1/2$.
This is consistent with the previous theoretical result \eqref{eq:spatial-error}.

\begin{figure}[htp]
\centering
      \includegraphics[width=4in,height=3in] {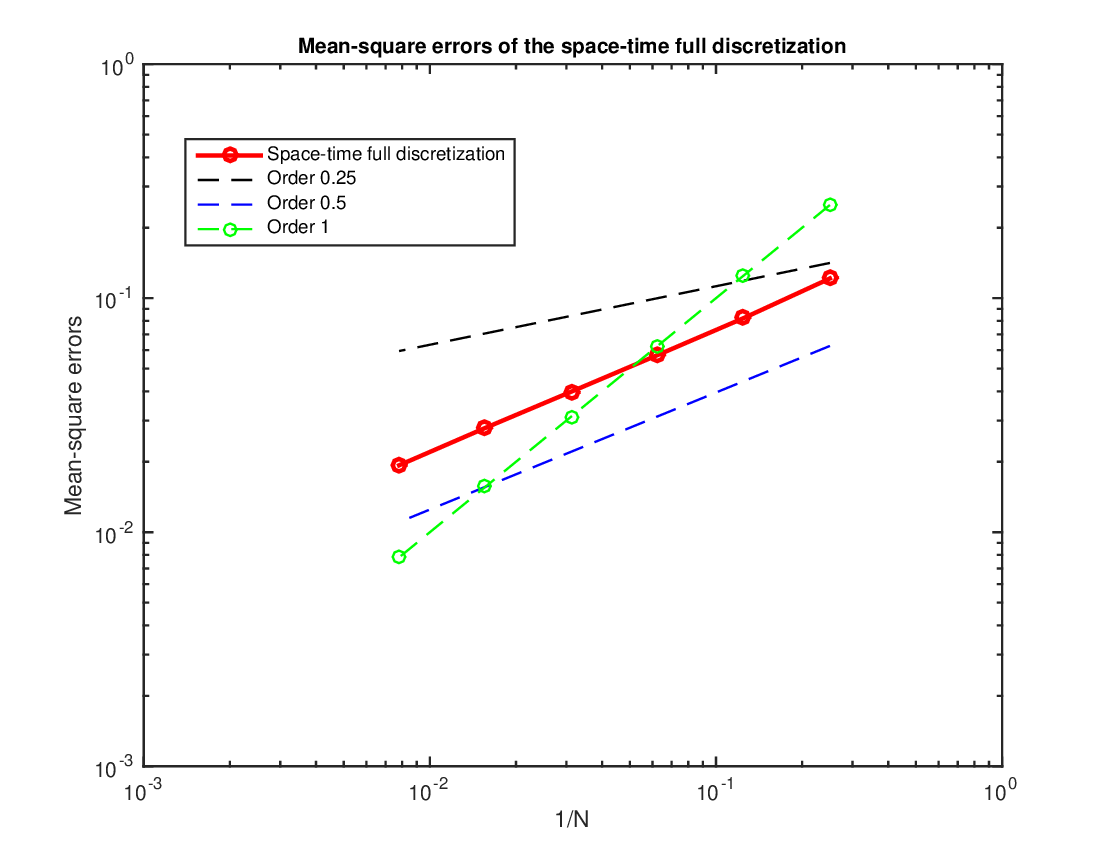}
 \caption{The convergence rate of the space-time full discretization.}
\label{fig:full-error}
\end{figure}

Moreover, we attempt to illustrate the error bound \eqref{eq:thm-main-result} for the fully discrete scheme \eqref{eq:num-result-AC-eqn}.
As implied by \eqref{eq:thm-main-result}, the convergence rate in space is identical to that in time. 
Consequently, we take $M = N$, $p = 2$  and $\beta = \frac12 - \epsilon$
with arbitrarily small $\epsilon >0$ in \eqref{eq:thm-main-result} to arrive at
\begin{equation} \label{eq:numerical-full-error}
 \| X(1) - Y^{N,N}_{t_N} \|_{L^2 ( \Omega; H) } \leq
C_{\epsilon} 
     N ^{ - \frac12 + \epsilon } .
\end{equation}
To see \eqref{eq:numerical-full-error}, we, similarly as above, do a full  discretization on a very fine mesh with
$M_{\text{exact}} = N_{\text{exact}} = 2^{11} = 2048$ to compute the ``exact'' solution.
Six different mesh parameters $N = 2^{i}, i = 2,3,...,7$ are then used to get six full discretizations. 
The resulting errors are listed in Table \ref{table:full-error} and plotted in Fig.\ref{fig:full-error} on a log-log scale. 
From Fig.\ref{fig:full-error},  one can observe the expected convergence rate of order almost $\tfrac12$,
which agrees with that indicated in \eqref{eq:numerical-full-error}.

%
%
%
%
   
\begin{table}[htp]
\begin{center} \footnotesize
\caption{Computational errors of the fully discrete scheme with $M =N$}
\label{table:full-error}
\begin{tabular*}{15cm}{@{\extracolsep{\fill}}cccccc}
\hline
$N = 2^2 $ & $ N = 2^3 $ & $N = 2^4$ & $N = 2^5$ & $N = 2^6$ & $N = 2^7$ \\
\hline
0.106381 & 0.077172 & 0.055174 & 0.039209  & 0.027624 & 0.019225 \\
\hline
\end{tabular*}
\end{center}
\end{table}

\bibliographystyle{abbrv}

\bibliography{bibfile}

\end{document}